\documentclass[11pt]{amsart}
\usepackage{txfonts}
\usepackage{} \textwidth=14.5cm \oddsidemargin=1cm
\usepackage{amsmath,yhmath}
\usepackage{amsxtra}
\usepackage{amscd}
\usepackage{amsthm}
\usepackage{amsfonts}
\usepackage{amssymb}
\usepackage{stmaryrd}
\usepackage{xcolor}

\input prepictex
\input pictex
\input postpictex

\newcommand{\newcom}{\newcommand}

\newcom{\al}{\alpha}
\newcom{\be}{\beta}
\newcom{\eps}{\epsilon}
\newcom{\del}{\delta}
\newcom{\Del}{\Delta}
\newcom{\ga}{\gamma}
\newcom{\Ga}{\Gamma}
\newcom{\ka}{\kappa}
\newcom{\Lam}{\Lambda}
\newcom{\lam}{\lambda}
\newcom{\Om}{\Omega}
\newcom{\om}{\omega}
\newcom{\Si}{\Sigma}
\newcom{\si}{\sigma}
\newcom{\tht}{\theta}
\newcom{\dtri}{\nabla}
\newcom{\td}{\tilde}
\newcom{\tri}{\triangle}
\newcom{\oo}{\infty}
\newcom{\vphi}{\varphi}
\newcom{\cB}{{\mathcal B}}
\newcom{\cC}{{\mathcal C}}
\newcom{\cD}{{\mathcal D}}
\newcom{\cF}{{\mathcal F}}

\newcom{\cL}{{\mathcal L}}

\newcom{\cK}{{\mathcal K}}
\newcom{\cP}{{\mathcal P}}
\newcom{\cR}{{\mathcal R}}
\newcom{\cS}{{\mathcal S}}
\newcom{\cU}{{\mathcal U}}
\newcom{\cX}{{\mathcal X}}
\newcom{\cY}{{\mathcal Y}}
\newcom{\cN}{{\mathcal N}}
\newcom{\cH}{{\mathcal H}}

\newcom{\R}{\Bbb R}
\newcom{\N}{\Bbb N}
\newcom{\Z}{\Bbb Z}
\newcom{\C}{\Bbb C}
\newcom{\E}{\Bbb E}

\newcom{\bx}{\bar x}
\newcom{\bz}{\bar z}
\newcom{\tx}{\tilde x}
\newcom{\tz}{\tilde z}
\newcom{\f}{\frac}
\newcom{\di}{\displaystyle\int}
\newcom{\ds}{\displaystyle\sum}
\newcom{\dl}{\displaystyle\lim}
\newcom{\ov}{\overline}
\newcom{\sset}{\subset}
\newcom{\wt}{\widetilde}
\newcom{\p}{\partial}
\newcom\na{\nabla}
\newcom{\co}{\cdot}
\newcom{\suml}{\sum\limits}
\newcom{\supl}{\sup\limits}
\newcom{\intl}{\int\limits}
\newcom{\infl}{\inf\limits}
\newcom{\disp}{\displaystyle}
\newcom{\non}{\nonumber}
\newcom{\no}{\noindent}
\newcom{\QED}{$\square$}
\def\ef{\hphantom{MM}\hfill\llap{$\square$}\goodbreak}

\newtheorem{athm}{\bf \t}[section]
\newenvironment{thm} [1] {\def\t{#1}\begin{athm} \bf \rm} {\end {athm}}
\newcom{\bthm}{\begin{thm}}\newcom{\ethm}{\end{thm}}

\newcom{\beq}{\begin{equation}}
\newcom{\eeq}{\end{equation}}
\newcom{\ben}{\begin{eqnarray}}
\newcom{\een}{\end{eqnarray}}
\newcom{\beno}{\begin{eqnarray*}}
\newcom{\eeno}{\end{eqnarray*}}

\numberwithin{equation}{section}

\begin{document}

\title[Weighted elliptic estimates]{Weighted elliptic estimates for a mixed boundary system related to the Dirichlet-Neumann operator on a corner domain}
\author{Mei Ming}
\address{School of Mathematics, Sun Yat-sen University, Guangzhou 510275, P. R. China}
\email{mingm@mail.sysu.edu.cn}

\date{Jan. 2019}
\maketitle


\begin{abstract}
Based on the $H^2$ existence of the solution, we investigate  weighted estimates for a mixed boundary elliptic system in a two-dimensional corner domain, when the contact angle $\om\in(0,\pi/2)$. This system is closely related to the Dirichlet-Neumann operator in the water-waves problem, and the weight we choose is decided by singularities  of the mixed boundary system.  Meanwhile, we also prove similar weighted estimates with a different weight for the Dirichlet boundary problem as well as the Neumann boundary problem when $\om\in(0,\pi)$.
\end{abstract}

\tableofcontents

\section{Introduction}

Based on the classical non-smooth theory in \cite{KMR, MR}, we consider   weighted estimates for a mixed boundary elliptic system in  a two-dimensional  corner domain $\Omega$. This domain is bounded by a top surface $\Ga_t=\{(x,z)\,|\,z=\eta(x)\}$  and a  smooth  bottom $\Ga_b=\{(x,z)\,|\,z=l(x)\}$, that is
\[\Omega=\{(x,z)\,|\,l(x)<z<\eta(x),\,x\ge 0\}\]  and 
the bottom satisfies
\[l(x)=-\ga x,\quad \hbox{when}\quad x\le x_0\] for some fixed constant $x_0$ and slope
$-\ga<0$.  Without loss of generality,  we place the contact point $X_c$ intersected by $\Ga_t$, $\Ga_b$ to be  at the origin $O=(0,0)$. The free surface  $z=\eta(x)$ satisfies
\[
\eta(0)=0,\quad\hbox{and}\quad 0<\eta(x)-l(x)\le H,\quad\forall x>0
\]
for some constant $H>0$.

Closely related to the Dirichler-Neumann operator in the water-waves problem,  we will focus on  the following mixed boundary elliptic problem for $u$ when proper conditions $h,\,f,\,g$ are given:
\[
\mbox{(MBVP)}\quad\left\{\begin{array}{ll}
\Delta u=h,\qquad \hbox{in}\quad \Omega\\
u\,|_{\Ga_t}=f,\qquad \p_{n_b}u\,|_{\Ga_b}=g.
\end{array}\right.
\]

When the domain changes with time $t$, the top surface becomes a free surface with a fixed bottom, and the contact point also varies. This kind of corner domains are related to a scene of sea waves moving near the beach in the real world, which are already used when we investigate the water-waves problem and related elliptic systems in \cite{MW1, MW2}.  For the moment, we only consider  a fixed  surface $z=\eta(x)$  independent of the time in this paper. 

\medskip
To prove the estimates for the mixed boundary  problem $\mbox{(MBVP)}$, one needs to notice firstly that this problem contains some singularity on the boundary, which requires naturally the non-smooth elliptic theory.  Therefore, before stating our main results, we shall recall some previous works on the non-smooth elliptic theory. 

 To start with,  here {\it non-smooth}  is generally referred to Lipschitz.  When the boundary is Lipschitz, the classical elliptic theory for a smooth boundary doesn't apply any more.
The non-smooth elliptic theory has  been fully developed in recent decades, and  fundamental works are done by Kondrat'ev \cite{Kon63, Kon67}. One can find some other early works by Birman and Skvortsov \cite{BS}, Eskin \cite{E}, Lopatinskiy \cite{Lop},  Maz'ya \cite{Maz64, Maz67}, Kondrat'ev and Oleinik \cite{KO}, Maz'ya and Plamenevskiy \cite{MazP77}, Maz'ya and Rossmann \cite{MazR92},  Grisvard \cite{PG1}, Dauge \cite{D} etc..
These works analyze singularities near the corner and provide regularity results in Sobolev space or weighted Sobolev space for general linear elliptic problems on Lipschitz domains. 

In fact, the existence of a variation solution in $H^1$  can be proved most of the time for a Lipschitz domain, see for example \cite{PG1}.  Compared to the smooth elliptic theory, when a higher regularity is considered,  the key for the non-smooth theory lies in  singularities, which can be expressed by a summation of singular functions  like $r^\lam\log^q r\varphi(\tht)$ near the corner point, where $r$ is the radius to the corner point, $\lam$ is an eigenvalue of the corresponding problem, $q$ is some constant, and $\varphi(\tht)$ is a bounded trigonometric function. Compared to  $H^1$ solutions,  it is well known that singularities arise when a higher-order regularity is referred to.  
At that time, the solution $u$ to an elliptic problem can be decomposed into  
\[
u=u_r+\sum_i c_{i}S_i
\] 
where $u_r$ is the regular part, $c_{i}$ the singular coefficient, and $S_i$ some singular function with an explicit formula as mentioned above. Moreover,  it is also well known that, the number in the summation of singular functions are finite and can be decided explicitly by the elliptic operator, the contact angle and the regularity, see for example \cite{PG1}. In fact, when one considers higher-order regularities or larger contact angles, the number  of singular functions  usually increases. 
The decompositions and  estimates for the regular part and the singular coefficients in Sobolev spaces  can be found in Kondrat'ev \cite{Kon67}, Maz'ya and Plamenevskiy \cite{MazP80, MazP84}, Dauge, Nicaise, Bourlard and Lubuma \cite{DNBL}, Grisvard \cite{PG1}, Costabel and Dauge \cite{CD93I, CD93II}, Ming and Wang \cite{MW1} etc.. 

\medskip
Based on the study for singular functions, a smart way to obtain a clean elliptic estimate  as in the classical case is to use weighted Sobolev spaces, for example, space $V^l_\be$ defined in Section 2 with some weight number $\be$  and order $l$. Due to the expressions of  singular functions, the weight is naturally in a form of $r^\be$, where  $r$ is the radius to the corner point. We refer to general weighted estimates in Kozlov, Mazya and Rossmann \cite{KMR}, Dauge \cite{D}, Mazya and Rossmann \cite{MR} etc.. These works provide some general weighted estimates assuming that the right side of the elliptic system also lies in corresponding weighted spaces. Meanwhile,  to obtain the weighted estimates, there are usually conditions between $\be, l$ and the eigenvalues of the corresponding eigenvalue problem: One requires that no eigenvalues $\lam$  lie  on the line $Re \lam=-\be+l-1$, see Theorem 6.1.1 \cite{KMR}.

\medskip

Using the weighted spaces introduced in these works and starting with the $H^2$ existence (which is already proved in \cite{MW1}), we prove {\it proper} weighted estimates for the mixed boundary problem and trace the dependence of the upper boundary in the coefficients.   

Firstly, ``{\it proper}" means that we identify the power $\be$ of the weight $r^\be$ very specifically, which is based on our analysis for the same mixed boundary problem in \cite{MW1}. On one hand,  the weight $\be$ we choose is decided by the order of singularities which appear in our problem. Thanks to Proposition 5.19 \cite{MW1}, when one considers $H^l(\Om)$  solution $u$ ($l\ge 3$), one needs  at least  the weight $r^{l-2}$ to eliminate the singular part 
\[
r^{-\f{(m+1/2)\pi}\om}\quad\hbox{ for}\  m\in \Z
\]
such that $r^{l-2}\na^l u\in L^2$ near the corner.
On the other hand, we obtain the weighted elliptic estimates without extra condition between eigenvalues $\lam$ and $\be, l$ as mentioned above (which is an important ingredient  in our results).
These two points result in the weighted space $V^l_{l-2+\be}(\Om)$ with $\be \in[0, 2]$ in our main theorem, which is defined in Section 2.

On the other hand,  one can see that the dependence of the upper boundary is not clearly proved in previous works. We provide detailed estimates for tracing this dependence  in this paper.

\medskip

The main theorem  is presented below.
\bthm{Theorem}\label{main thm}(Mixed boundary)
{\it Let the contact angle $\om\in(0,\pi/2)$ and $u\in H^2(\Om)$ be the solution to $\mbox{(MBVP)}$ for given $h\in L^2(\Om)$, $f\in H^{3/2}(\Ga_t)$ and $g\in H^{1/2}(\Ga_b)$.  Moreover, for a real $\be\in [0,2]$ and  an integer $l\ge 2$, one assumes  that 
\[
h\in V^{l-2}_{l-2+\be}(\Om),\quad f\in V^{l-1/2}_{l-2+\be}(\Ga_t)\quad\hbox{and}\quad g\in V^{l-3/2}_{l-2+\be}(\Ga_b).
\]
Then, one has $u\in V^l_{l-2+\be}(\Om)$, and \\
\noindent (i) If $\eta\in W^{l,\infty}(\R^+)$,  there holds
\[
\|u\|_{V^l_{l-2+\be}(\Om)}\le C(\|\eta'\|_{W^{l-1,\infty}(\R^+)})\big(\|h\|_{V^{l-2}_{l-2+\be}(\Om)}+\|f\|_{V^{l-1/2}_{l-2+\be}(\Ga_t)}+\|g\|_{V^{l-3/2}_{l-2+\be}(\Ga_b)}\big);
\]
\noindent (ii) If $\eta\in H^{l-1/2}(\R^+)$ and $l\ge 3$, there holds
\[
\|u\|_{V^l_{l-2+\be}(\Om)}\le C(\|\eta\|_{H^{l-1/2}(\R^+)})\big(\|h\|_{V^{l-2}_{l-2+\be}(\Om)}+\|f\|_{V^{l-1/2}_{l-2+\be}(\Ga_t)}+\|g\|_{V^{l-3/2}_{l-2+\be}(\Ga_b)}\big).
\]
The coefficient $C$ is a positive polynomial of $\|\eta'\|_{W^{l-1,\infty}(\R^+)}$ or $\|\eta\|_{H^{l-1/2}(\R^+)}$.
}
\ethm

\bthm{Remark}\label{rmk}{\it In fact, this result can be adjusted immediately to the case of a bounded corner domain, where there are two contact points between the upper surface and the bottom.
}
\ethm

On the other hand, we  also consider about weighted estimates for the Dirichlet boundary problem 
\[
\mbox{(DVP)}\quad\left\{\begin{array}{ll}
\Delta u=h,\qquad \hbox{in}\quad \Omega\\
u\,|_{\Ga_t}=f,\qquad u\,|_{\Ga_b}=g
\end{array}\right.
\]
as well as  for the Neumann boundary problem
\[
\mbox{(NVP)}\quad\left\{\begin{array}{ll}
\Delta u=h,\qquad \hbox{in}\quad \Omega\\
\p_{n_t}u\,|_{\Ga_t}=f,\qquad \p_{n_b}u\,|_{\Ga_b}=g
\end{array}\right.
\]
with the compatibility condition 
\[
\int_\Om h dX=\int_{\Ga_t}fds+\int_{\Ga_b}g ds.
\]

Similar weighted estimates are proved in this paper for both $\mbox{(DVP)}$ and $\mbox{(NVP)}$ when the contact angle varies in a much larger interval, and meanwhile the weight is slightly different from the mix-boundary case. 
\bthm{Theorem}\label{D and N thm} {\it Assume that the contact angle $\om\in(0,\pi)$,  $\be\in (0,1]$ be a real number and the integer $l\ge 2$.\\
(i) (Dirichlet boundary)   Let $u\in H^2(\Om)$ be the solution to $\mbox{(DVP)}$ for given $h\in L^2(\Om)$, $f\in H^{3/2}(\Ga_t)$ and $g\in H^{3/2}(\Ga_b)$ satisfying
\[
f|_{X_c}=g|_{X_c}.
\]  Moreover, one assumes  that 
\[
h\in V^{l-2}_{l-1+\be}(\Om),\quad f\in V^{l-1/2}_{l-1+\be}(\Ga_t)\quad\hbox{and}\quad g\in V^{l-1/2}_{l-1+\be}(\Ga_b).
\]
Then, one has $u\in V^l_{l-1+\be}(\Om)$.  When $\eta\in H^{l-1/2}(\R^+)$ and $l\ge 3$, there holds
\[
\|u\|_{V^l_{l-1+\be}(\Om)}\le C(\|\eta\|_{H^{l-1/2}(\R^+)})\Big(\|h\|_{V^{l-2}_{l-1+\be}(\Om)}+\|f\|_{V^{l-1/2}_{l-1+\be}(\Ga_t)}+\|g\|_{V^{l-1/2}_{l-1+\be}(\Ga_b)}\Big).
\]
(ii) (Neumann boundary)  Let $u\in H^2(\Om)$ be the solution to $\mbox{(NVP)}$ for given $h\in L^2(\Om)$, $f\in H^{1/2}(\Ga_t)$ and $g\in H^{1/2}(\Ga_b)$.  Moreover, one assumes  that 
\[
h\in V^{l-2}_{l-1+\be}(\Om),\quad f\in V^{l-3/2}_{l-1+\be}(\Ga_t)\quad\hbox{and}\quad g\in V^{l-3/2}_{l-1+\be}(\Ga_b).
\]
Then, one has $u\in V^l_{l-1+\be}(\Om)$.  When $\eta\in H^{l-1/2}(\R^+)$ and $l\ge 3$, there holds
\[
\|u\|_{V^l_{l-1+\be}(\Om)}\le C(\|\eta\|_{H^{l-1/2}(\R^+)})\Big(\|h\|_{V^{l-2}_{l-1+\be}(\Om)}+\|f\|_{V^{l-3/2}_{l-1+\be}(\Ga_t)}+\|g\|_{V^{l-3/2}_{l-1+\be}(\Ga_b)}\Big).
\]
The coefficient $C$ above is a positive polynomial of $\|\eta\|_{H^{l-1/2}(\R^+)}$.
}
\ethm

\subsection{Organization of the paper} 
In Section 2 we introduce the weighted spaces on $\Om$ and its boundaries with some useful lemmas. Section 3 proves the main theorem for the mixed boundary problem. In Section 4, some other boundary problems are considered, while Section 5 provides the application of our theory on the Dirichlet-Neumann operator.

\subsection{Notations}

\noindent - $X_c$ denotes the contact point. We simply set $X_c=O(0,0)$ here;\\
\noindent - $\cK$ is the cone $\{(x,z)\,|\,-\ga x\le z\le\eta'(0)x\}$;\\
\noindent - $\cC$ is the strip $\{(t,\tht)\,|\,t\in \R,\,-\om_2\le \tht\le\om_1\}$;\\
\noindent - The angular interval $I=[-\om_2,\om_1]$. The contact angle  $\om=\om_1+\om_2$;\\
\noindent - $\Ga_t,\Ga_b$ denote the upper boundary and the lower boundary respectively for the domain $\Om$, $\cK$ or $\cC$, when no confusion will be made;\\
\noindent - Recalling from \cite{MW1} that, the function $d(\cdot)$ introduced in the transformation of the domain is 
\[
d(\cdot)=1-\f 1{\eta'(\bar\eta^{-1}(\cdot))+\ga}
\] 
with $\bar \eta( x)=\eta(x)+\ga  x$ invertible near $X_c$.

\bigskip

\section{Weighted Sobolev spaces on corner domains}

\subsection{Definitions for weighted spaces and transformations of domains} 
We will introduce definitions of weighted Sobolev spaces firstly on the cone $\cK$ and then on the corner domain $\Om$, which can be found in \cite{KMR, MR}. 

For an integer $l\ge 0$ and a real $\be$, the space $V^l_\be(\cK)$ can be defined as the closure of  $C^\infty_0(\bar \cK \backslash X_c)$ with the norm
\[
\|w\|_{V^l_\be(\cK)}=\Big(\displaystyle \int_\cK \sum_{|\al|\le l} r^{2(\be-l+|\al|)}\big|\na^\al_X\,w\big|^2dX\Big)^\f12
\]
with $r$ the radius with respect to $X_c$.

Next, we recall  straightening transformations $T_S$ and $T_R$ from \cite{MW1}.
To begin with, Let 
\[
\cS=\{(\tilde x,\tilde z)\,|\,\tilde x\ge 0,\,0\le \tilde z\le \tilde x\}.
\]
 $T_S$ is the local transformation near the point $X_c$ which maps  $\cS\cap U_{\del \cS}$  into  $\Om\cap U_\del$:
\[T_S: \quad (\tilde x,\,\tilde z)\in \cS\cap U_{\del \cS}\mapsto (\bar x,\,\bar z)\in\Om\cap U_\del\]
with
\[\bar x=\tilde x+\bar\eta^{-1}(\tilde z)- \tilde z,\quad \bar z=\tilde z-\ga\Big(\tilde x+\bar\eta^{-1}(\tilde z)- \tilde z\Big)\]
where $\bar\eta^{-1}(\tilde z)$ is the inverse of $\bar\eta(\bar x)=\eta(\bar x)+\ga  \bar x$. $U_{\del \cS}$ and $U_{\del}$ are two corresponding neighborhoods of $X_c$. We know that $T_S$ is invertible:
\[T_S^{-1}:\quad (\bar x,\,\bar z)\in \Om\cap U_\del\mapsto (\tilde x,\,\tilde z)\in \cS\cap U_{\del \cS}\] 
where
\[\tilde x=\bar x-\bar\eta^{-1}(\ga\bar x+\bar z)+\ga \bar x+\bar z,\quad \tilde z=\ga\bar x+\bar z.\] 
Moreover, we introduce the linear transform 
\[
T_0:\quad X=(x,\,z)\in \cK\mapsto \tilde X=(\tilde x,\,\tilde z)= X P_0\in\cS
\] 
with 
\[
P_0=\left(\begin{matrix} 1+\ga \,d(0)& \ga\\ d(0) & 1\end{matrix}\right)\quad\hbox{where}\quad d(0)=1-\f{1}{\ga+\eta'(0)}.
\] Together with $T_S$, we set
\[
T_c=T_S\circ T_0
\] which maps the cone $\cK$ to the domain $\Om$ near the corner.

Besides, we also have the transform $T_R$ which maps a flat strip $R$ to the rest part of $\Om$:
\[T_R:\quad (x,\,z)\in R\mapsto (\bar x,\,\bar z)\in \Om\]
with
\[\bar x= x,\quad \bar z=\eta(x)z+l(x)(1-z),\] where $R=\{(x,\,z)|\, x\ge x_\del,\,0\le z\le 1\}$ is a flat strip, and $x_\del>0$ is a constant fixed by  $U_\del$.  The inverse transform $T_R$ is
\[
T_R^{-1}:\quad (\bar x,\,\bar z)\in \Om\cap \{\bar x\ge x_\del\}\mapsto (x,\,z)\in R
\] 
where
\[x=\bar x,\quad z=\f{\bar z-l(\bar x)}{\eta(\bar x)-l(\bar x)}.\]

Now it's the time to define the weighted space $V^l_\be(\Om)$ on  $\Om$.  We firstly set $\chi_c\in C^\infty_0(\bar\Om)$ supported near $X_c$ with some diameter $\del>0$ small enough. Since singularities only take place near the corner point $X_c$, the weight also concentrates near the corner. 
The weighted space $V^l_\be(\Om)$ is equipped with the norm 
\beq\label{def of V}
\|u\|_{V^l_\be(\Om)}=\|v_c\|_{V^l_\be(\cK)}+\|v_R\|_{H^l(R)}
\eeq where
\[
v_c=u_c\circ T_c\quad\hbox{with}\quad u_c=\chi_c u
\] and 
\[
v_R=(1-\chi_c)u\circ T_R.
\]
Obviously, the space doesn't depend on the choices of the cut-off function $\chi_c$.

On the other hand, one also needs to use another type of weighted space $W^l_{2,\be}(\cC)$ on the infinite strip $\cC=\R\times[-\om_2,\,\om_1]$, which can be found in \cite{KMR}. In fact, for a function  $w(t,\tht)$ on $\cC$, the norm for $W^l_{2,\be}(\cC)$ is defined as 
\beq\label{W norm}
\|w\|_{W^l_{2,\be}(\cC)}=\|e^{\be t}w\|_{H^l(\cC)}.
\eeq
Similarly, the corresponding weighted space $W^{l-1/2}_{2,\be}(\Ga_t)$, $W^{1-3/2}_{2,\be}(\Ga_b)$ on the upper and lower boundaries are defined with norms
\beq\label{W boundary norm}
\|w\|_{W^{l-1/2}_{2,\be}(\Ga_t)}=\|e^{\be t}w\|_{H^{l-1/2}(\Ga_t)},\quad \|w\|_{W^{l-3/2}_{2,\be}(\Ga_b)}=\|e^{\be t}w\|_{H^{l-3/2}(\Ga_b)}.
\eeq
Moreover, $W^l_{2,\be}(\R)$ used in Section 2.3 is defined in a similar way. 

In the end, we recall a regularizing diffeomorphism near $X_c$ from \cite{MW1} which is a variation based on the transformations $T_S$ and $T_c$.  

To begin with, we define $\tilde s(\tilde x,\tilde z)$  on $\cS$ satisfying the Dirichlet boundary condition:
\[
\tilde s(\tilde x,\tilde z)\big|_{\Ga_t:\,\tilde z=\tilde x}=\be(\tilde x)\bar\eta^{-1}(\tilde x) 
\] 
where $\be$ is a cut-off function defined on $[0,+\infty)$ and vanish away from $0$.  
Consequently, one has from Remark 4.8 \cite{MW1} that if $\be\bar \eta^{-1}\in H^{l-1/2}(\R^+)$, then $\tilde s(x,z)\in H^{l}(\cS)$ with the estimate
\beq\label{trace for s}
\|\tilde s\|_{H^l(\cS)}\le C \,|\be \bar \eta^{-1}|_{H^{l-1/2}(\R^+)}\le C\big(|\eta|_{H^{m_0}(\R^+)}\big) |\eta|_{H^{l-1/2}(\R^+)},
\eeq 
where the constant $m_0> 3/2$.

As a result, we define the regularized transformation $\tilde T_S$ as 
\[\tilde T_S: \quad (\tilde x,\,\tilde z)\in S\cap U_{\del S}\mapsto (\bar x,\,\bar z)\in\Om\cap U_\del\]
with
\[\bar x=\tx+\tilde s\big(\eps \tx+(1-\eps) \tz,\,\tz\big)- \tz,\quad \bar z=\tz-\ga\Big(\tx+\tilde s \big(\eps \tx+(1-\eps)\tz,\,\tz\big)-\tz\Big)\] where   $\eps$ is a small constant to be explained.
A direct computation shows that
\[
Det(\na \tilde T_s)=1+\eps\, \p_{\tx} \tilde s,
\]   so $\tilde T_S$ is invertible as long as  the constant $\eps$ is  small enough such that 
\[\eps\le \f 1{2\|\p_{\tx}\tilde s \|_\infty}.\]
 
Some more computations lead to the associated coefficient matrix related to $\mbox{(MBVP)}$ 
\[
\tilde P_S=(\na \tilde T^{-1}_S)\circ \tilde T_S=\f 1{1+\eps\,\p_{\tx} \tilde s}\left(\begin{matrix}
1+\ga\big(1-(1-\eps)\p_{\tx} \tilde s-\p_{\tz} \tilde s\big) & \ga(1+\eps\,\p_{\tx}\tilde s)\\
1-(1-\eps)\p_{\tx} \tilde s-\p_{\tz} \tilde s & 1+\eps\,\p_{\tx} \tilde s
\end{matrix}\right)\Big|_{(\eps \tx+(1-\eps) \tz,\,\tz)} 
\] 
and we denote 
\[
\tilde P_0=\tilde P_S|_{X_c}=\left(\begin{matrix} a&b\\c &d\end{matrix}\right).
\]
Similarly as before, the transformation $\tilde T_0$ from $\cK$ to $\cS$ is defined as 
\[
\tilde T_0:\quad X=(x,z)\in\cK\mapsto \tilde X=(\tx,\tz)=XP_0\in\cS,
\]
and we also define on $\cK$ that
\beq\label{def of s}
s(x,z)=\tilde s(XP_0P_\eps )=\tilde s(\eps\tx+(1-\eps)\tz,\tz)\quad\hbox{with}\quad 
P_\eps=\left(\begin{matrix} \eps&0\\1-\eps &1\end{matrix}\right).
\eeq
So one can replace $P_S$ in system \eqref{vc system} by
\beq\label{tilde PS}
\tilde P_S\circ \tilde T_0=\f 1{1+a\p_xs+c\p_zs}\left(\begin{matrix}
1+\ga\big(1-b\p_xs-d\p_zs\big) & \ga(1+a\p_xs+c\p_zs)\\
1-b\p_xs-d\p_zs & 1+a\p_xs+c\p_zs
\end{matrix}\right)
\eeq
when we need it. 

Similarly as before, we set
\beq\label{def of tilde Tc}
\tilde T_c=\tilde T_S\circ \tilde T_0
\eeq
which maps the cone $\cK$ to the domain $\Om$ near the corner.

\subsection{Traces on the boundary}
The weighted spaces on the boundary are also needed in our theory. We introduce the definitions of the trace spaces  (see \cite{MR}), and some trace theorems are  discussed, too.

Firstly, we define $V^{l-1/2}_\be(\Ga_t)$ (and $V^{l-1/2}_\be(\Ga_b)$) for $l\ge 1$ as the spaces for traces of functions from $V^l_\be(\Om)$ on $\Ga_t$ ( and $\Ga_b$) respectively. 

$V^{l-1/2}_\be(\Ga_t)$ is equipped with the norm
\[
\|u\|_{V^{l-1/2}_\be(\Ga_t)}=\inf \big\{\|u_{ex}\|_{V^l_\be(\Om)}\,\big|\,u_{ex}\in V^l_\be(\Om),\,u_{ex}|_{\Ga_t}=u\big\},
\]
and the norm of $V^{l-1/2}_\be(\Ga_b)$ is defined similarly. 

Consuquently, one concludes the following lemma immediately.
\bthm{Lemma}\label{restriction trace}
{\it Let $u\in V^l_\be(\Om)$ and set $f=u|_{\Ga_t}$, then one has $f\in V^{l-1/2}_\be(\Ga_t)$ and the estimate
\[
 \|f\|_{V^{l-1/2}_\be(\Ga_t)}\le \|u\|_{V^l_\be(\Om)}.
\]
 Similar conclusion holds for the trace on $\Ga_b$ and for the case $V^l_\be(\cK)$.
}
\ethm

\medskip
Since the traces related to the cone $\cK$ will be used frequently, one needs to go further with the norms defined above. Notice that the angle $\theta\equiv\om_1$ on the upper boundary $\Ga_t$ of $\cK$,  and $\theta\equiv-\om_2$ for $\Ga_b$.  Lemma 6.1.2 from \cite{KMR} gives an equivalent norm for $V^{l-1/2}_\be(\Ga_t)$ as 
\beq\label{boundary norm}
\begin{split}
\|u\|^2_{V^{l-1/2}_\be(\Ga_t)}=&\sum_{j\le l-1}\displaystyle\int_{\R^+}r^{2(\be-l)+1}\big|(r\p_r)^j u(r,\om_1)\big|^2dr\\
&\quad
+\sum_{j\le l-1}\displaystyle\int_{\R^+}\displaystyle\int_{\R^+}r^{2(\be-l)+2}\f{\big|(r\p_r)^{j}u(r,\om_1)-(\rho\p_\rho)^{j}u(\rho,\om_1)\big|^2}{|r-\rho|^2}d\rho dr,
\end{split}
\eeq
which will be used frequently in our paper. The equivalent norm for $V^{l-1/2}_\be(\Ga_b)$ is  defined similarly.

The following lemma concerns the trace theorem with Dirichlet boundary conditions, which is modified from Lemma 2.2.1 \cite{MR}.
\bthm{Lemma}\label{Dirichlet trace thm}(Dirichlet boundary) {\it Let $f\in V^{l-1/2}_\be(\Ga_t)$ and $g\in V^{l-1/2}_\be(\Ga_b)$ with integer $l\ge 1$. Then there exists a function $w\in V^{l}_\be(\cK)$ such that 
\[
w|_{\Ga_t}=f,\quad w|_{\Ga_b}=g
\]
with the estimate
\[
\|w\|_{V^l_\be(\cK)}\le C\big(\|f\|_{V^{l-1/2}_\be(\Ga_t)}+\|g\|_{V^{l-1/2}_\be(\Ga_b)}\big)
\]
where the constant $C$ depends only  on $\be,l,\cK$.
}
\ethm
\noindent{\bf Proof}. Taking $m=1$ in Lemma 2.2.1 \cite{MR}, one obtains the desired result immediately.
\ef

\subsection{Some premilinaries }
Some preparations are done in this part. Firstly, embeddings between different weighted spaces are discussed. Moreover, one considers the relationships between different weighted spaces and ordinary spaces. In the end, the Laplace transform is  introduced with some basic properties, and an equivalent norm for a weighted space is defined based on this transform.

The functions considered here are always compactly supported near $X_c$ with a size $\del$, and we focus on the cone $\cK$ most of the time.
\bthm{Lemma}\label{imedding}
{\it 
Assume that integers $l_2\ge l_1\ge 0$ and real $\be_1,\be_2$ satisfy
\[
l_2-\be_2\ge l_1-\be_1.
\]
For any $v\in V^{l_2}_{\be_2}(\cK)$ with a compact support of size $\del$ near $X_c$, one can have $v\in V^{l_1}_{\be_1}(\cK)$  such that
\[
\|v\|_{V^{l_1}_{\be_1}(\cK)}\le \del^{(l_2-\be_2)-(l_1-\be_1)}\|v\|_{V^{l_2}_{\be_2}(\cK)}.
\]
 Moreover, similar results hold for $V^{l-1/2}_{\be}(\Ga_t)$ and $V^{l-3/2}_{\be}(\Ga_b)$ with constants $C=C(l_1,l_2,\be_1,\be_2,\del)$.
}
\ethm
\noindent{\it Proof}. 
One only needs to check from the definitions to prove this lemma.  In fact, for any $v\in V^{l_2}_{\be_2}(\cK)$ with a compact support of size $\del$ near $X_c$, a simple computation  shows that
\[
\begin{split}
\|r^{(\be_1-l_1)+|\al|}\p^\al v\|_{L^2(\cK)}&=\|r^{(l_2-\be_2)-(l_1-\be_1)}r^{(\be_2-l_2)+|\al|}\p^\al v\|_{L^2(\cK)}\\
&\le \del^{(l_2-\be_2)-(l_1-\be_1)}\|r^{(\be_2-l_2)+|\al|}\p^\al v\|_{L^2(\cK)}
\end{split}
\]
where $\p^\al=\p_x^{\al_1}\p_z^{\al_2}$ satisfying $|\al|=\al_1+\al_2\le l_1$.  Therefore, the case for 
$V^{l_2}_{\be_2}(\cK)$ is proved, and the other cases can be done similarly.
\ef

\bthm{Lemma}\label{H2 V22}
{\it
Let $v$ and $f$ be two functions on $\cK$ and $\Ga_t$  (or $\Ga_b$) respectively with a compact support of size $\del$ near $X_c$. \\
\noindent(i) When $v\in H^2(\cK)$, one has $v\in V^2_2(\cK)$ satisfying
\[
\|v\|_{V^2_2(\cK)}\le C\|v\|_{H^2(\cK)};
\]
\noindent(ii) When $v\in L^2(\cK)$, one has $v\in V^0_2(\cK)$ satisfying 
\[
\|v\|_{V^0_2(\cK)}\le C\|v\|_{L^2(\cK)};
\]
\noindent(iii) When $f\in H^{3/2}(\Ga_t)$, one has $f\in V^{3/2}_2(\Ga_t)$ satisfying 
\[
\|f\|_{V^{3/2}_2(\Ga_t)}\le C\|f\|_{H^{3/2}(\Ga_t)}.
\] 
Similar inequality holds also for the case from $H^{1/2}(\Ga_b)$ to $V^{1/2}_2(\Ga_b)$.
Moreover, the constant $C$ above depends on $\del,\cK$.
}
\ethm
\noindent{\it Proof}. The first two cases can be proved in a similar way as in Lemma \ref{imedding}, and it only remains to check (iii).

In fact, when $f\in V^{3/2}_2(\Ga_t)$, one knows directly from the definition that
\[
\begin{split}
\|f\|^2_{V^{3/2}_2(\Ga_t)}&=\displaystyle\int_{\R^+}r|f(r)|^2dr+\displaystyle\int_{\R^+}r^3|f'(r)|^2dr+\displaystyle\int_{\R^+}\int_{\R^+}r^2\f{|f(r)-f(\rho)|^2}{|r-\rho|^2}drd\rho\\
&\qquad+\displaystyle\int_{\R^+}\int_{\R^+}r^2\f{|rf'(r)-\rho f'(\rho)|^2}{|r-\rho|^2}drd\rho\\
&\triangleq A_1+A_2+A_3+A_4
\end{split}
\]
where the first two terms can be handled easily since $f$ is compactly supported near $X_c$:
\[
A_1+A_2\le  (\del+\del^3)\|f\|^2_{H^1(\Ga_t)}.
\]
Now it remains to take care of the last two terms.  To begin with, one has 
\[\begin{split}
A_3
&\le \displaystyle\int_{\R^+}\int^{2\rho}_{\rho/2}\dots drd\rho+\displaystyle\int_{\R^+}\int^{\rho/2}_{0}\dots drd\rho+\displaystyle\int_{\R^+}\int^{+\infty}_{2\rho}\dots drd\rho\\
&\triangleq A_{31}+A_{32}+A_{33},
\end{split}\]
where a direct analysis shows that
\[
A_{31}\le \del^2\, C\displaystyle\int_{\R^+}\int^{2\rho}_{\rho/2}\f{|f(r)-f(\rho)|^2}{|r-\rho|^2}drd\rho\le \del^2\, C\|f\|^2_{H^{1/2}(\Ga_t)}\]
since one has $r\sim \rho$ in this case and remember that $f$ is compactly supported near $X_c$. Moreover, one can also  have
\[\begin{split}
A_{32}&\le C\Big(\displaystyle\int_{\R^+}\int^{\infty}_{2r}r^2|f(r)|^2\f{1}{|r-\rho|^2}d\rho\,dr+\displaystyle\int_{\R^+}\int^{\rho/2}_0|f(\rho)|^2\f{r^2}{|r-\rho|^2}dr\,d\rho\Big)\\
&\le C\displaystyle\int_{\R^+} r|f(r)|^2dr\le \del \,C\|f\|^2_{H^1(\Ga_t)},
\end{split}\]
and a similar inequality holds for $A_{33}$. Consequently, we arrive at 
\[
A_3\le \del \,C\|f\|^2_{H^{3/2}(\Ga_t)}.
\]
On the other hand, similar computations can be done for the term $A_4$. Therefore,  the proof for the case $ H^{3/2}(\Ga_t)$ is finished.
\ef

\bthm{Lemma}\label{H32 V32_0}
{\it Let $f\in V^{3/2}_0(\Ga_t)$ and $g\in V^{1/2}_0(\Ga_b)$ be functions compactly supported near $X_c$ of $\cK$ with a size $\del$. Then one has $f\in H^{3/2}(\Ga_t)$ and $g\in H^{1/2}(\Ga_b)$ satisfying 
\[
\|f\|_{H^{3/2}(\Ga_t)}\le C\|f\|_{V^{3/2}_0(\Ga_t)},\quad \|g\|_{H^{1/2}(\Ga_b)}\le C\|g\|_{V^{1/2}_0(\Ga_b)}
\]
where the constant $C$ depends on $\del,\cK$.
}\ethm
\noindent{\it Proof}. The proof can be done similarly as in the previous lemma.  In fact,  using the definition of $H^{3/2}(\Ga_t)$ and $V^{3/2}_0(\Ga_t)$, one writes directly that
\[
\|f\|^2_{H^{3/2}(\Ga_t)}=\|f\|^2_{H^1(\Ga_t)}+\int_{\R^+}\int_{\R^+}\f{|f'(r)-f'(\rho)|^2}{|r-\rho|^2}drd\rho
\]
and 
\[
\begin{split}
\|f\|^2_{V^{3/2}_0(\Ga_t)}=&\|r^{-3/2}f\|^2_{L^2(\Ga_t)}+\|r^{-1/2}f'\|^2_{L^2(\Ga_t)}+\int_{\R^+}\int_{\R^+}r^{-2}|\f{|f(r)-f(\rho)|^2}{|r-\rho|^2}drd\rho\\
&\quad +\int_{\R^+}\int_{\R^+}r^{-2}\f{|rf'(r)-\rho f'(\rho)|^2}{|r-\rho|^2}drd\rho.
\end{split}
\]
Since $f$ is supported near $X_c$ with a size $\del$, one can easily see that
\[
\|f\|^2_{H^1(\Ga_t)}=\|r^{3/2}r^{-3/2} f\|^2_{L^2(\Ga_t)}+\|r^{1/2}r^{-1/2}f'\|^2_{L^2(\Ga_t)}\le \del C\|f\|^2_{V^{3/2}_0(\Ga_t)},
\]
so it remains to check the last term in $H^{3/2}$ norm. 

Similarly as before, one can write 
\[
\begin{split}
\displaystyle\int_{\R^+}\int_{\R^+}\f{|f'(r)-f'(\rho)|^2}{|r-\rho|^2}drd\rho&=\displaystyle\int_{\R^+}\int^{2\rho}_{\rho/2}\dots drd\rho+\displaystyle\int_{\R^+}\int^{\rho/2}_{0}\dots drd\rho+\displaystyle\int_{\R^+}\int^{+\infty}_{2\rho}\dots drd\rho\\
&\triangleq A_1+A_2+A_3.
\end{split}
\]
Direct computations show that
\[
\begin{split}
A_1&=\displaystyle\int_{\R^+}\int^{2\rho}_{\rho/2}r^{-2}\f{|rf'(r)-rf'(\rho)|^2}{|r-\rho|^2}drd\rho\\
&\le C\Big(\displaystyle\int_{\R^+}\int^{2\rho}_{\rho/2}r^{-2}\f{|rf'(r)-\rho f'(\rho)|^2}{|r-\rho|^2}drd\rho+\displaystyle\int_{\R^+}\int^{2\rho}_{\rho/2}r^{-2}|f'(\rho)|^2drd\rho\Big)
\le C\|f\|^2_{V^{3/2}_0(\Ga_t)},
\end{split}
\]
and $A_2,A_3$ can be handled similarly as before. Consequently, the case of $H^{3/2}(\Ga_t)$ is proved. Moreover, the case of $H^{1/2}(\Ga_b)$ can also be proved similarly.
\ef

The following lemma deals with the relationship between $V^l_\be(\cK)$ and $W^l_{2,\be}(\cC)$, which is quoted directly from (6.1.6) and (6.1.7)  \cite{KMR}.
\bthm{Lemma}\label{space V and W}
{\it Let $r=e^t$ with $(r, \tht)$ polar coordinates  and denote $w(t,\tht)=v(r,\tht)$, where $v(r,\tht)$ is defined on $\cK$. Then $w(t,\tht)$ is defined on $\cC$ and there exist constants $C_1,\,C_2$ depending on $l,\be$ and $\cK$ such that
\[
C_1 \|w\|_{W^l_{2,\be-l+1}(\cC)}\le \|v\|_{V^l_\be(\cK)}\le C_2\|w\|_{W^l_{2,\be-l+1}(\cC)}\quad\hbox{i.e.}\quad \|v\|_{V^l_\be(\cK)}\simeq\|w\|_{W^l_{2,\be-l+1}(\cC)}.
\]
Moreover, similar results hold on the boundary:
\[
\|v\|_{V^{l-1/2}_\be(\Ga_t)}\simeq\|w\|_{W^{l-1/2}_{2,\be-l+1}(\Ga_t)}\quad \hbox{and}\quad \|v\|_{V^{l-3/2}_\be(\Ga_b)}\simeq\|w\|_{W^{l-3/2}_{2,\be-l+1}(\Ga_b)}.
\]
}\ethm

In the end of this section, we introduce the Laplace transform $\cL$ acting on any $w\in C^\infty_0(\R)$:
\[
\breve w(\lam)=\cL w(\lam)=\displaystyle\int_{\R}e^{-\lam t}w(t)dt, \quad \forall \lam\in \C.
\]
Some well-known properties of this transform are recalled below, quoted directly from Lemma 5.2.3 \cite{KMR}.
\bthm{Lemma}\label{laplace transform}{\it (i) The transform defines a linear and continuous mapping from $C^\infty_0(\cR)$ into the space of analytic functions on the complex plane $\C$. Further more, one has
\[
\cL(\p_t w)=\lam \cL w
\]
\noindent (ii) For all $u,v\in C^\infty_0(\R)$, the Parseval equality 
\[
\displaystyle\int^{+\infty}_{-\infty}e^{2\be t}u(t)\,\overline{v(t)}dt=\f{1}{2\pi i}\displaystyle\int_{Re \lam=-\be}\breve u(\lam)\,\overline{\breve v(\lam)}d\lam
\]
holds. Here the integration on the right side takes place over the line $l_{-\be}=\{\lam=i\tau-\be,\ \tau\in\R\}$, and $\bar z$ denotes the conjugation of $z$.\\
\noindent (iii) The inverse Laplace transform is given by the formula
\[
w(t)=(\cL^{-1}\breve w)(t)=\f{1}{2\pi i}\displaystyle\int_{Re \lam=-\be} e^{\lam t}\breve w(\lam)d\lam.
\]
\noindent (iv) If $w\in W^0_{2,\be_1}(\R)\cap W^0_{2,\be_2}(\R)$ where $\be_1<\be_2$, then $\breve w=\cL w$ is holomorphic on the strip $-\be_2<Re \lam <-\be_1$.
}\ethm

Combining these properties with the definition of $W^l_{2,\be}(\cC)$ and $W^{l-1/2}_{2,\be}(\Gamma_t)$,  we recall from \cite{KMR} the following lemma.
\bthm{Lemma}\label{equivalent W norm}
{\it The norm \eqref{W norm} with an integer $l\ge 0$ is equivalent to the norm 
\[
\|w\|=\Big( \f{1}{2\pi i}\displaystyle\int_{Re\lam=-\be}\|\breve w\|^2_{H^l(I,\lam)}d\lam\Big)^{1/2}
\]
where 
\[
\|\breve w\|_{H^l(I,\lam)}=\Big(\|\breve w(\lam,\cdot)\|^2_{H^l(I)}+|\lam|^{2l}\|\breve w(\lam,\cdot)\|^2_{L_2(I)}\Big)^{1/2}
\]
Analogously, an equivalent norm to \eqref{W boundary norm} for $W^{l-1/2}_{2, \be}(\Ga_t)$ ($l\ge1 $) on $\Ga_t$ is 
\[
\|w\|=\Big( \f{1}{2\pi i}\displaystyle\int_{Re\lam=-\be}(1+|\lam|^{2l-1})|\,\breve w(\lam,\om_1)|^2d\lam\Big)^{1/2},
\]
and there is a similar norm for the case of $W^{l-3/2}_{2, \be}(\Ga_b)$.
}\ethm

\section{Estimates for the mixed boundary problem}
We start with the existence of the solution to  $\mbox{(MBVP)}$ in certain weighted space, and then the regularity is considered. The weighted estimates is proved by an induction argument, and the dependence of the upper boundary is traced at the same time.

To begin with, one must consider about the existence of the solution in proper weighted space, which we wish to be built on the existence result in  ordinary Sobolev spaces from \cite{BR, MW1}.
In fact, recalling Theorem 5.2 and Remark 5.3 \cite{MW1}, we state the following lemma for the unique existence of the solution in $H^2(\Om)$.
\bthm{Lemma}\label{H2 existence}
{\it Suppose that $\Om$ has a $\cC^{2,0}$ upper boundary $\Ga_t$. Let functions $h\in L^2(\Om)$, $f\in H^{3/2}(\Ga_t)$ and $g\in H^{1/2}(\Ga_b)$ be given and the contact angle $\om\in(0,\pi/2)$. Then there exists a unique solution $u\in H^2(\Om)$ to $\mbox{(MBVP)}$.
}
\ethm

One can see from the definition of $V^l_\be(\Om)$ that,  the parts of the norm near the corner and away from the corner are treated in completely different ways. We use weighted spaces near the corner, and the elliptic estimates  need to be proved (which is the key ingredient in our paper). When it is away from the corner, ordinary Sobolev spaces are used, so standard elliptic estimates can be applied directly. In a word, we focus on the weighted estimates near the corner  in the following text. 

Recalling from \eqref{def of V} that,  we defined on $\cK$ the function
\[
v_c=u_c\circ T_c\qquad\hbox{with}\quad u_c=\chi_c u
\]
where  $u\in H^2(\Om)$ the solution to $\mbox{(MBVP)}$.

Some computations as in \cite{MW1} show that $v_c$ satisfies the system
\beq\label{vc system}
\left\{\begin{array}{ll}\na\cdot P_c\na v_c=h_c\qquad \hbox{on}\quad \cK\\
v_c|_{\Ga_t}= f_c,\quad \p^{P_c}_{n_b}v_c|_{\Ga_b}= (1+\ga^2)^{1/2}g_c\end{array}\right.
\eeq where $-\ga$ is the constant slope of $\Ga_b$ near $X_c$ and 
\[
h_c=\big(\chi_c \,h-[\chi_c,\Del]u\big)\circ T_c,\quad f_c=\big(\chi_c|_{\Ga_t}\,f\big)\circ T_c,\quad g_c=\big(\chi_c|_{\Ga_b}\,g-(\p_{n_b}\chi_c)u|_{\Ga_b}\big)\circ T_c.
\]
Besides,  the coefficient matrix is
\[
P_c=(P^{-1}_0)^tP^t_SP_SP^{-1}_0
\] 
where
\[
P_S=(\na T^{-1}_S)\circ T_S=\left(\begin{matrix} 1+\ga \,d(\ga x+z)& \ga\\ d(\ga x+z) & 1\end{matrix}\right),\quad
P_0=P_S\big|_{X_c}=\left(\begin{matrix} 1+\ga \,d(0)& \ga\\ d(0) & 1\end{matrix}\right).
\]
\medskip

To prove the main theorem, we need to focus on system \eqref{vc system} for $v_c$. First of all, under the assumptions of Theorem \ref{main thm} and combining Lemma \ref{H2 existence}, one finds immediately  that 
\beq\label{v H2 regularity}
v_c\in H^2(\cK)\quad\hbox{with}\quad h_c\in L^2(\cK),\ f_c\in H^{3/2}(\Ga_t),\ g_c\in H^{1/2}(\Ga_b)
\eeq
while notice that all functions are compactly supported near $X_c$.

Now we are in a position to  introduce proper weighted spaces for the system of $v_c$. In fact, combining Lemma \ref{H2 V22}, one has immediately 
\beq\label{basic regularity}
v_c\in V^2_2(\cK),\  h_c\in V^0_2(\cK),\  f_c\in V^{3/2}_2(\Ga_t),\ \hbox{and}\ g_c\in V^{1/2}_2(\Ga_b)
\eeq
where the weight $\be=2$. 
Based on these spaces, we  improve the regularity of $v_c$ in the following two subsections. One will see that, when the contact angle $\om\in(0,\pi/2)$ and a proper weight $\be$ is chosen, there is no extra singularity when higher regularity is considered.

\subsection{Lower-order regularity near the corner}

So far,  $v_c$ belongs to $V^2_2(\cK)$ with the weight $\be=2$. The aim of this subsection is to show that $v_c$ also belongs to $V^2_0(\cK)$ with a lower weight $\be=0$, which is a very important step and leads us to the proper weighted space $V^l_{l-2}(\cK)$. 
\bthm{Proposition}\label{weighted estimate near Xc}
{\it Let $v_c$ be the solution to \eqref{vc system} and \eqref{v H2 regularity} holds. Moreover,  for a real $\be\in [0,2]$ one assumes  that 
\[
h_c\in V^{0}_{\be}(\cK),\quad f_c\in V^{3/2}_{\be}(\Ga_t),\quad g_c\in V^{1/2}_{\be}(\Ga_b),
\]
and 
\[
\|\eta\|_{W^{2,\infty}}\le C_0
\]
for some constant $C_0$.

Then one has $v_c\in V^2_{\be}(\cK)$ and the weighted estimate holds
\beq\label{vc V20 estimate}
\|v_c\|_{v^2_{\be}(\cK)}\le C\big(\|h_c\|_{V^{0}_{\be}(\cK)}+\|f_c\|_{V^{3/2}_{\be}(\Ga_t)}+\|g_c\|_{V^{1/2}_{\be}(\Ga_b)}\big),
\eeq
where the constant $C=C(\cK,\be)$.
}
\ethm

Before we prove this proposition, some preparations are needed.  Firstly, let 
\[
\cB=\big(\na\cdot P_c\na, \ \cdot|_{\Ga_t},\ \p^{P_c}_{n_b}\cdot|_{\Ga_b}\big)
\]
be the elliptic operator for  system \eqref{vc system}. In particular, since direct computations show that 
\[
\cP_c|_{X_c}=Id,
\]
  we denote
\[
\cB_0=\cB\big|_{X_c}=\big(\Del,\ \cdot|_{\Ga_t},\ \p_{n_b}\cdot|_{\Ga_b}\big)
\] 
as the restriction of coefficients in $\cB$ on the contact point $X_c$.

System \eqref{vc system} can be rewritten into a perturbation form of operator $\cB_0$ near the contact point $X_c$:
\[
\cB_0\, v_c=(h_c,\,f_c,\,g_c)-(\cB-\cB_0)v_c
\]
or equivalently the following system
\beq\label{vc system perturbation}
\left\{\begin{array}{ll}\Del v_c=h_v\qquad \hbox{on}\quad \cK\\
v_c|_{\Ga_t}= f_v,\quad \p_{n_b}v_c|_{\Ga_b}= g_v\end{array}\right.
\eeq 
where 
\[
h_v=h_c-\na\cdot(P_c-Id)\na v_c,\ f_v=f_c,\ \hbox{and}\ g_v=(1+\ga^2)^{1/2}g_c-\p^{P_c-Id}_{n_b}v_c\big|_{\Ga_b}.
\]
\bthm{Lemma}\label{vc system right side}
{\it Under the assumptions of Proposition \ref{weighted estimate near Xc} with $\be\in [0,2)$, there exists  a number $\eps\in(0,1)$ depending on $\be$ such that
\[
h_v\in V^0_{1+\eps}(\cK),\ f_v\in V^{3/2}_{1+\eps}(\cK),\ \hbox{and}\ g_v\in V^{1/2}_{1+\eps}(\Ga_b)
\]
in system \eqref{vc system perturbation}.
}
\ethm
\noindent{\it Proof}. Since one assumes that $h_c\in V^0_\be(\cK)$,   applying Lemma \ref{imedding} with $l_2=l_1=0$, $\be_2=\be$ and $\be_1=1+\eps$ on $h_c$   leads to 
\[
h_c\in V^0_{1+\eps}(\cK).
\]
Here one requires that $\be \le 1+\eps$.
Similarly,  applying Lemma \ref{imedding} with $l_2=l_1=2$, $\be_2=\be$ and $\be_1=1+\eps$ on $f_c$ and $g_c$ leads to
\[
f_c\in V^{3/2}_{1+\eps}(\Ga_t),\ g_c\in V^{1/2}_{1+\eps}(\Ga_b).
\]
It remains to deal with perturbation terms 
\[
(\cB-\cB_0)v_c=\big(\na\cdot (P_c-Id)\na v_c,\ 0,\ \p^{P_c-Id}_{n_b}v_c|_{\Ga_b}\big).
\]
Firstly,  one can show directly that 
\[
r^{1+\eps}\,\na\cdot(P_c-Id)\na v_c\in L^2(\cK)
\]
since one has $v_c\in H^2(\cK)$ with a compact support near $X_c$ and the assumption for $\eta$ in  Proposition \ref{weighted estimate near Xc}. 
This infers that $\na\cdot(P_c-Id)\na v_c\in V^0_{1+\eps}(\cK)$.

On the other hand, the boundary term can be written as 
\[
\p^{P_c-Id}_{n_b}v_c|_{\Ga_b}={\bf n}_b\cdot (P_c-Id)\na v_c\big|_{\Ga_b},
\]
so one can show in a similar way as above that
\[
(P_c-Id)\na v_c\in V^1_{1+\eps}(\cK).
\]
Applying Lemma \ref{restriction trace} on $\Ga_b$, one has immediately that $\p^{P_c-Id}_{n_b}v_c|_{\Ga_b}\in V^{1/2}_{1+\eps}(\Ga_b)$. 

Summing up these results above, one can finish the proof.
\ef

\medskip
For the moment,  we are ready to change the weight for $v_c$, that is, from $V^2_2(\cK)$ to $V^2_{1+\eps}(\cK)$. Concerning elliptic systems on corner domains, it is well known  that one will meet with singularities most of the time when one wants to consider about two different spaces, see for example \cite{PG, KMR}. The key lemma below tells us that in our settings with the contact angle $\om\in (0,\pi/2)$, no singularity happens when we choose the space carefully. Moreover, we only investigate about the proper space without establishing any estimate at this time.
\bthm{Proposition}\label{singularity decomposition}
{\it Let the contact angle $\om\in (0,\pi/2)$. Assume that system \eqref{vc system perturbation} admits a solution $v_c\in V^2_2(\cK)$ with 
\[
h_v\in V^0_{1+\eps}(\cK)\cap V^0_2(\cK),\quad f_v\in V^{3/2}_{1+\eps}(\Ga_t)\cap V^{3/2}_2(\Ga_t)\quad \hbox{and}\ g_v\in V^{1/2}_{1+\eps}(\Ga_b)\cap V^{1/2}_2(\Ga_b),
\]
then one has $v_c\in V^2_{1+\eps}(\cK)$ without any singularity decomposition.
}
\ethm
\noindent{\it Proof}. 
The idea of this proof  follows the proofs for Theorem 5.4.1 and Theorem 6.1.4 \cite{KMR}. In fact, we convert system \eqref{vc system} on the cone $\cK$ equivalently to a system on a horizontal strip, and  then the Laplace transform is  applied to derive  the related eigenvalue problem. As a result, the solution $v_c$ under Laplace transform could be expressed through an ODE. Based on some analysis on  eigenvalues, we are able to use Cauchy's Formula to show that $v_c$ eventually lies in the desired weighted space.\\
\noindent{\bf Step 1}. Change of variable.  First of all, system \eqref{vc system perturbation} can be rewritten under polar coordinates:
\[
\left\{\begin{array}{ll}(r^2\p^2_r+\p^2_\tht+r\p_r) v_c=r^2\,h_v\qquad \hbox{on}\quad \cK\\
v_c|_{\tht=\om_1}= f_v,\quad -\p_\tht v_c|_{\tht=-\om_2}=  r\,g_v.\end{array}\right.
\]
Secondly, introducing the following change of variable 
\[
t=\ln r,\quad\hbox{i.e.}\quad r=e^t\quad\hbox{for} \quad \forall\,t\in \R
\] and denoting
\[
w(t,\tht)=v_c(r,\tht),
\]
 the system above for $v_c$ can be changed equivalently into the system for $w(t,\tht)$ on an infinite strip $\cC=\R\times[-\om_2,\,\om_1]$:
\beq\label{w system}
\left\{\begin{array}{ll}(\p^2_t+\p^2_\tht) w=e^{2t}\,h_w\qquad \hbox{on}\quad \cC\\
w|_{\tht=\om_1}= f_w,\quad -\p_\tht w|_{\tht=-\om_2}=  e^t\,g_w\end{array}\right.
\eeq 
with the notations 
\[
h_w(t,\tht)=h_v(r,\tht),\ f_w(t,\tht)=f_v(r,\tht),\ \hbox{and}\ g_w(t,\tht)=g_v(r,\tht).
\]
As a result,   applying the assumptions of this proposition and  Lemma \ref{space V and W} on $v_c$, $h_v$, $f_v$ and $g_v$, one derives immediately that
\[
w\in W^2_{2,1}(\cC),
\] and the right side of \eqref{w system} satisfies
\beq\label{right side of w system}
e^{2t}h_w\in W^0_{2,\eps}(\cC)\cap W^0_{2,1}(\cC),\ f_w\in W^{3/2}_{2,\eps}(\Ga_t)\cap W^{3/2}_{2,1}(\Ga_t),\  e^t\,g_w\in W^{1/2}_{2,\eps}(\Ga_b)\cap W^{1/2}_{2,1}(\Ga_b).
\eeq
\noindent{\bf Step 2}. Laplace transform to an ODE. 
One performs the Laplace transform on $w(t,\tht)$ with respect to $t$ and denote
\[
\breve w(\lam,\cdot)=(\cL w)(\lam,\cdot)=\displaystyle\int_{\R}e^{-\lam t}w(t,\cdot)dt, \quad \forall \lam\in \C
\]
Applying Lemma \ref{laplace transform} (i), one arrives at the system for $\breve w(\lam,\tht)$ from \eqref{w system}:
\beq\label{bw system}
\left\{\begin{array}{ll}\lam^2\breve w+\p^2_\tht \breve w=\cL(e^{2t}\,h_w),\qquad \tht\in I\\
\breve w|_{\tht=\om_1}= \cL(f_w),\quad -\p_\tht \breve w|_{\tht=-\om_2}=\cL\big( e^t\,g_w),\end{array}\right.
\eeq 
and one knows from \eqref{right side of w system} and Lemma \ref{equivalent W norm} that 
\beq\label{L(h,f,g) space}
\breve w\in H^2(I,\lam),\quad\cL(e^{2t}\,h_w)\in L_2(I,\lam).
\eeq
One can see that our system \eqref{vc system} turns into an ordinary differential system with parameter $\lam$, which becomes more handy. 

We denote by 
\[
\cU(\lam)=\big(-\p^2_\tht-\lam^2, \,\cdot|_{\tht=\om_1},\,-\p_\tht\cdot|_{\tht=-\om_2}\big)
\]
  the operator  of  system \eqref{bw system} with parameter $\lam\in\C$. For each fixed $\lam$, $\cU(\lam)$ continuously maps $H^l(I)$ into $H^{l-2}(I)$ for any $l\ge 2$ with corresponding boundary values.

A direct computation shows that the corresponding eigenvalue problem for $\cU(\lam)$ reads
\[
\left\{\begin{array}{ll}
-\phi''(\tht)-\lam^2\phi(\tht)=0,\quad \tht\in I\\
\phi(\om_1)=0,\qquad -\phi'(-\om_2)=0
\end{array}\right.
\]
where the eigenvalues are  countable and real with the explicit expressions
\[
\lam_m=\f{(m+1/2)\pi}{\om}\qquad\hbox{for}\ \forall m\in\Z.
\]
By the way, the eigenfunctions are $\phi_m(\tht)=\cos{\big(\lam_m(\tht+\om_2)\big)}$. In fact, these eigenvalues and eigenfunctions coincide with those in \cite{MW1}, which is characteristic for the mixed-type elliptic problem.

Since  the contact angle $\om$ is assumed to be in $(0,\pi/2)$ in this paper, one finds immediately that
\[
\lam_m\notin [-1,1]\quad \forall m\in\Z,
\]
which implies 
\[
\cU(\lam)\ \hbox{ is invertible when }\  \lam\in[-1,1].
\]

\noindent{\bf Step 3}. Singularity decomposition without singularity.
For this moment, we plan to show that $w\in W^2_{2,0}(\cC)$ by solving system \eqref{bw system}.
First of all, we will start from expressing $\breve w\in H^2(I,\lam)$ in terms of the right hand side.  

In fact, when we take $Re\lam=-1$, system \eqref{bw system} is uniquely solvable in $ H^2(I)$ (which is already known since system \eqref{w system} admits a solution $w\in W^2_{2,1}(\cC)$ ). Moreover,  one can express the solution as below 
\[
\breve w(\lam,\tht)=\cU(\lam)^{-1}\big(\cL(e^{2t}\,h_w),\,\cL(f_w),\,\cL\big( e^t\,g_w)\big)
\in H^2(I,\lam),
\]
where  $\cL(e^{2t}\,h_w)$ satisfies \eqref{L(h,f,g) space} with $Re\lam=-1$. 

Applying the inverse Laplace transform  and Lemma \ref{laplace transform} (iii), one obtains
\beq\label{w expression}
w(t,\tht)=\f{1}{2\pi i}\displaystyle\int_{Re\lam=-1}e^{\lam t}\,\breve w(\lam,\tht)d\lam\in W^2_{2,1}(\cC).
\eeq
From Lemma \ref{laplace transform} (iv), one can see that for each $\tht\in I$,  $\cL(e^{2t}\,h_w)$, $\cL(f_w)$ and $\cL\big( e^t\,g_w)$ are holomorphic in the strip $-1<Re\lam<-\eps$. 
Therefore, the only singularities of the function 
\[
e^{\lam t}\,\breve w(\lam,\tht)=e^{\lam t}\,\cU(\lam)^{-1}\big(\cL(e^{2t}\,h_w),\,\cL(f_w),\,\cL\big( e^t\,g_w)\big)
\]
from \eqref{w expression} in the strip $-1<Re\lam<-\eps$ are the poles of $\cU(\lam)^{-1}$, i.e. the eigenvalues of $\cU(\lam)$. Combining previous analysis on $\cU(\lam)$, this implies immediately that no singularity takes place in the strip $-1<Re\lam<-\eps$.
\medskip

Now we are in a position to show that $w\in W^2_{2,0}(\cC)$. In fact,  let $\rho>0$ to be a constant, then the complex  domain
\[
D_\rho=\big\{\lam\in\C\,|\, -1<Re\lam <-\eps, \ |Im\lam|>\rho\big\}
\]
doesn't contain any eigenvalue of $\cU(\lam)$. 

Rewriting \eqref{w expression} and  applying Cauchy's Formula, we have
\[
\begin{split}
w(t,\tht)=&\f{1}{2\pi i}\lim_{\rho\rightarrow +\infty}\displaystyle\int^{-1+i\rho}_{-1-i\rho}e^{\lam t}\,\breve w(\lam,\tht)d\lam\\
=& \f{1}{2\pi i}\lim_{\rho\rightarrow +\infty}\Big(\displaystyle\int^{-\eps+i\rho}_{-\eps-i\rho}e^{\lam t}\,\breve w(\lam,\tht)d\lam+\displaystyle\int^{-\eps-i\rho}_{-1-i\rho}e^{\lam t}\,\breve w(\lam,\tht)d\lam-\displaystyle\int^{-\eps+i\rho}_{-1+i\rho}e^{\lam t}\,\breve w(\lam,\tht)d\lam\Big).
\end{split}
\]
A lemma is needed here to deal with the last two integrals, which will be proved after the proof of this proposition. 
\bthm{Lemma}\label{limit zero}{\it 
Let $\cU(\lam)$ be the operator  of system \eqref{bw system} for $\breve w(\lam,\tht)$ with $(h_w,f_w,g_w)$ satisfying \eqref{right side of w system}. Then there holds
\[
\lim_{\rho\rightarrow +\infty}\displaystyle\int^{-\eps\pm i\rho}_{-1\pm i\rho}e^{\lam t}\,\cU(\lam)^{-1}\big(\cL(e^{2t}\,h_w),\,\cL(f_w),\,\cL\big( e^t\,g_w)\big)d\lam=0.
\]
}
\ethm
 Consequently, with the help of this lemma we arrive at
\[
w(t,\tht)=\f{1}{2\pi i}\lim_{\rho\rightarrow +\infty}\displaystyle\int^{-\eps+i\rho}_{-\eps-i\rho}e^{\lam t}\,\breve w(\lam,\tht)d\lam
=\f{1}{2\pi i}\displaystyle\int_{Re\lam=-\eps}e^{\lam t}\,\breve w(\lam,\tht)d\lam.
\]
Recalling that $\breve w\in H^2(I,\lam)$, we can finally conclude with Lemma  \ref{laplace transform} (ii) (iii) that
\[
w\in W^2_{2,\eps}(\cC).
\]
Therefore,  we  apply Lemma \ref{space V and W} to find 
\[
v_c\in V^2_{1+\eps}(\cK)\cap V^2_2(\cK)
\]
and the proof is finished.
\ef

\noindent{\bf Proof of Lemma \ref{limit zero}}.  Since this proof is adapted from the proof of Lemma 5.4.1 \cite{KMR},  we only sketch the main idea here to be self-content.  

Firstly, let
\[
w_\rho(t,\tht)=\displaystyle\int^{-\eps+i\rho}_{-1+i\rho}e^{\lam t}\,\cU^{-1}(\lam)\big(\cL(e^{2t}\,h_w),\,\cL(f_w),\,\cL\big( e^t\,g_w)\big)d\lam=\displaystyle\int^{-\eps+i\rho}_{-1+i\rho}e^{\lam t}\,\breve w(\lam,\tht)d\lam.
\]
Taking $L^2$ norm on $\cC_N=[-N,N]\times I $ for a constant $N>0$, one has 
\[
\|w_\rho\|^2_{L^2(\cC_N)}
= \displaystyle\int_{\cC_N}|w_\rho(t,\tht)|^2dtd\tht
\le C\displaystyle\int_I\displaystyle\int^{-\eps+i\rho}_{-1+i\rho}\big|\breve w(\lam,\tht)\big|^2d\lam d\tht
\]
On the other hand, checking from Theorem 3.6.1 \cite{KMR} one finds the elliptic estimate for system \eqref{bw system} of $\breve w(\lam,\tht)$:
\[
\|\breve w\|_{H^l(I,\lam)}\le C\Big(\|\cL(e^{2t}\,h_w)\|_{H^{l-2}(I,\lam)}+(1+|\lam|^{l-1/2})|\cL(f_w)|+(1+|\lam|^{l-3/2})|\cL\big( e^t\,g_w)|\Big),
\]
which implies immediately 
\[
\begin{split}
&\displaystyle\int_I\big|\breve w(\lam,\tht)\big|^2d\tht=\|\breve w(\lam,\cdot)\|^2_{L^2(I)}\\
&\le C|\lam|^{-2l}\Big(\|\cL(e^{2t}\,h_w)\|^2_{H^{l-2}(I,\lam)}+(1+|\lam|^{2l-1})|\cL(f_w)|^2+(1+|\lam|^{2l-3})|\cL\big( e^t\,g_w)|^2\Big).
\end{split}
\]
As a result, one can show that
\[
\begin{split}
&\displaystyle\int^{c_2}_{c_1}\|w_\rho\|^2_{L^2(\cC_N)}d\rho\\
& \le C \displaystyle\int^{c_2}_{c_1}\displaystyle\int^{-\eps+i\rho}_{-1+i\rho}\Big(\|\cL(e^{2t}\,h_w)\|^2_{H^{l-2}(I,\lam)}+(1+|\lam|^{2l-1})|\cL(f_w)|^2+(1+|\lam|^{2l-3})|\cL\big( e^t\,g_w)|^2\Big)d\lam d\rho
\end{split}
\]
with the constant $C=C(N,c_1)$. Rewriting this double integral by changing the order of the integration, one derives 
\[
\begin{split}
&\displaystyle\int^{c_2}_{c_1}\|w_\rho\|^2_{L^2(\cC_N)}d\rho\\
&  \le C \displaystyle\int^{1}_{\eps}\displaystyle\int_{Re\lam=-\be}\Big(\|\cL(e^{2t}\,h_w)\|^2_{H^{l-2}(I,\lam)}+(1+|\lam|^{2l-1})|\cL(f_w)|^2+(1+|\lam|^{2l-3})|\cL\big( e^t\,g_w)|^2\Big)d\lam d\be,
\end{split}
\]
which together with Lemma \ref{equivalent W norm} leads to 
\[
\displaystyle\int^{c_2}_{c_1}\|w_\rho\|^2_{L^2(\cC_N)}d\rho \le C \displaystyle\int^{1}_{\eps}\Big(\|e^{2t}h_w\|^2_{W^{l-2}_{2,\be}(\cC)}+\|f_w\|^2_{W^{l-1/2}_{2,\be}(\Ga_t)}+\|e^tg_w\|^2_{W^{l-3/2}_{2,\be}(\Ga_b)}\Big)d\be.
\]
Therefore, combining \eqref{right side of w system}, one knows that $\|w_\rho\|_{L^2(\cC_N)}$ is also square integrable over the interval $(c_1,\infty)$ and the proof can be finished.
\ef

\medskip

In order to prove the estimate in Proposition \ref{weighted estimate near Xc},  we quote the weighted elliptic estimate for system \eqref{vc system perturbation} in the following lemma,  which can be found in Theorem 6.1.1 \cite{KMR}. Notice that this lemma  holds due to the previous analysis on $\cU(\lam)$: No eigenvalues of $\cU(\lam)$  lie on the line $Re\lam =-\be+l-1$, where we take $\be\in [0,2]$ and $l=2$ here.

\bthm{Lemma}\label{Del system estimate}
{\it Let $\be\in[0,2]$. Assume that  there exists a solution $v_c\in V^2_\be(\cK)$ for system \eqref{vc system perturbation}  with $h_v\in V^{0}_\be(\cK)$, $f_v\in V^{3/2}_{\be}(\Ga_t)$ and $g_v\in V^{1/2}_\be(\Ga_b)$.  Then there holds
\[
\|v_c\|_{V^2_\be(\cK)}\le C\big(\|h_v\|_{V^{0}_\be(\cK)}+\|f_v\|_{V^{3/2}_{\be}(\Ga_t)}+\|g_v\|_{V^{1/2}_\be(\Ga_b)}\big),
\]
where the constant $C=C(\cK,\be )$.
}
\ethm

Now we are ready to prove Proposition \ref{weighted estimate near Xc}.\\
\noindent{\bf Proof for Proposition \ref{weighted estimate near Xc}}. 
Firstly, one needs to show that $v_c\in V^2_\be(\cK)$. The case when $\be=2$ has been proved in \eqref{basic regularity}. For the case when $\be\in [0,2)$,   applying Lemma \ref{vc system right side}, one can see that the assumptions of Proposition \ref{singularity decomposition} are satisfied  for some $\eps\in (0,1)$  depending on $\be$. As a result, one knows from Proposition \ref{singularity decomposition} that  $v_c\in V^2_{1+\eps}(\cK)$.  Repeating this procedure finite times to reach a lower weight at each time, one can finally show that $v_c\in V^2_0(\cK)$. 

Secondly, to prove the weighted estimate \eqref{vc V20 estimate}, one applies Lemma \ref{Del system estimate} on  system \eqref{vc system perturbation} to derive 
\[
\|v_c\|_{V^2_\be(\cK)}\le C\big(\|h_v\|_{V^0_\be(\cK)}+\|f_v\|_{V^{3/2}_\be(\Ga_t)}+\|g_v\|_{V^{1/2}_\be(\Ga_t)}\big)
\]
with the constant $C=C(\cK,\be)$.  Substituting the expressions of $h_v,f_v,g_v$ from system \eqref{vc system perturbation}, one has 
\beq\label{vc V20 estimate 1}
\begin{split}
\|v_c\|_{V^2_\be(\cK)}\le & C\big(\|h_c\|_{V^0_\be(\cK)}+\|f_c\|_{V^{3/2}_\be(\Ga_t)}+\|g_c\|_{V^{1/2}_\be(\Ga_t)}\\
&\qquad +\|\na\cdot (P_c-Id)\na v_c\|_{V^0_\be(\cK)}+\|\p^{P_c-Id}_{n_b}v_c\|_{V^{1/2}_\be(\Ga_b)}\big),
\end{split}
\eeq
where the last two terms need to be handled. 

In fact, recall from system \eqref{vc system} that $P_c=(P^{-1}_0)^t P^t_S P_S P^{-1}_0$,  which implies 
\[
P_c-Id=(P^{-1}_0)^t(P_S-P_0)^tP_SP^{-1}_0+(P^{-1}_0)^tP^t_0(P_s-P_0)^tP^{-1}_0
\]
where
\[
P_S-P_0=\big(d(z)-d(0)\big)\left(\begin{matrix} \ga & 0\\ 1 & 0\end{matrix}\right)\quad\hbox{with}\quad
d(z)-d(0)=\f{\eta'\big(\bar\eta^{-1}(z)\big)-\eta'(0)}{(\eta'(0)+\ga)\big(\eta'(\bar\eta^{-1}(z))+\ga\big)}.
\] Here, recall from \cite{MW1} that we have $\bar\eta^{-1}(0)=0$ since we set $\eta(0)=0$.

Consequently, one can show directly that 
\[
\begin{split}
\|\na\cdot (P_c-Id)\na v_c\|_{V^0_\be(\cK)}&=\|r^\be\na\cdot (P_c-Id)\na v_c\|_{L^2(\cK)}\\
&\le \del \,C(\|\eta\|_{W^{2,\infty}})\|r^\be\na^2 v_c\|_{L^2(\cK)}+C(\|\eta\|_{W^2,\infty})\|r^\be\na v_c\|_{L^2(\cK)}\\
&\le \del\, C(\|\eta\|_{W^{2,\infty}})\|v_c\|_{V^2_\be(\cK)}
\end{split}
\]
where $\del$ comes from $d(z)-d(0)$ and remember that $v_c$ is compactly supported near $X_c$ with radius $\del$. Moreover, the last step is proved using the inequality
\[
\|r^\be\na v_c\|_{L^2(\cK)}=\|r\, r^{\be-1}\na v_c\|_{L^2(\cK)}\le \del \|v_c\|_{V^2_\be(\cK)}.
\]
Secondly, one has for the term $\|\p^{P_c-Id}_{n_b}v_c\|_{V^{1/2}_\be(\Ga_b)}$ the following estimate
\[
\begin{split}
\|\p^{P_c-Id}_{n_b}v_c\|_{V^{1/2}_\be(\Ga_b)}\le & C\|\p^{P_c-Id}_{n_b}v_c\|_{V^1_\be(\cK)}\\
\le & C\big(\|r^{\be-1} \p^{P_c-Id}_{n_b}v_c\|_{L^2(\cK)}+\|r^\be\na \p^{P_c-Id}_{n_b}v_c\|_{L^2(\cK)}\big)
\end{split}
\]
where the constant vector ${\bf n}_b$ is extended on $\cK$ and  Lemma \ref{restriction trace} is applied. Besides, one can show similarly as before that 
\[
\begin{split}
\|r^{\be-1} \p^{P_c-Id}_{n_b}v_c\|_{L^2(\cK)}&=\|r^{\be-1} {\bf n}_b\cdot (P_c-Id)\na v_c\|_{L^2(\cK)}\\
&\le \del \,C(\|\eta\|_{W^2,\infty}) \|v_c\|_{V^2_\be(\cK)}
\end{split}
\]
and 
\[
\|r^\be\na \p^{P_c-Id}_{n_b}v_c\|_{L^2(\cK)}\le \del \,C(\|\eta\|_{W^2,\infty}) \|v_c\|_{V^2_\be(\cK)}.
\]
Summing these up, one arrives at 
\[
\|\p^{P_c-Id}_{n_b}v_c\|_{V^{1/2}_\be(\Ga_b)}\le \del \,C(\|\eta\|_{W^2,\infty}) \|v_c\|_{V^2_\be(\cK)}.
\]
As a result, substituting the estimates above into \eqref{vc V20 estimate 1}, we conclude  that 
\[
\|v_c\|_{V^2_\be(\cK)}\le  C\Big(\|h_c\|_{V^0_\be(\cK)}+\|f_c\|_{V^{3/2}_\be(\Ga_t)}+\|g_c\|_{V^{1/2}_\be(\Ga_t)} +\del \,C(\|\eta\|_{W^2,\infty}) \|v_c\|_{V^2_\be(\cK)}\Big).
\]
In the end, using the assumption that $\|\eta\|_{W^2,\infty}\le C_0$, the proof can be finished if one choose $\del$ small enough depending on $C_0$.
\ef

\subsection{higher-order regularity near the corner}
In this part, we continue to improve the regularity of $v_c$. At this time, the target space is $V^l_{l-2+\be}(\cK)$ for integer $l\ge 2$. The case when $l=2$ is already considered in last subsection. Compared to previous analysis, we don't meet with singularity here, and standard elliptic theory can be applied locally for the regularity.
\bthm{Proposition}\label{order-l regularity}
{\it Let   $l\ge3$, $\be\in [0,2]$ and the contact angle $\om\in (0,\pi/2)$. Assume that system \eqref{vc system perturbation} admits a solution $v_c\in V^{l-1}_{l-3+\be}(\cK)$ with 
\[
h_c\in V^{l-2}_{l-2+\be}(\cK),\ f_c\in V^{l-1/2}_{l-2+\be}(\Ga_t)\  \hbox{and}\ g_c\in V^{l-3/2}_{l-2+\be}(\Ga_b),
\]
then one has $v_c\in V^{l}_{l-2+\be}(\cK)$ with the estimate
\[
\|v_c\|_{V^{l}_{l-2+\be}(\cK)}\le C(\|\eta'\|_{W^{l-1,\infty}})\Big(\|h_c\|_{V^{l-2}_{l-2+\be}(\cK)}+\|f_c\|_{V^{l-1/2}_{l-2+\be}(\Ga_t)}+\|g_c\|_{V^{l-3/2}_{l-2+\be}(\Ga_b)}\Big).
\]
}
\ethm
\noindent{Proof}. Firstly, we will prove the case when $l=3$. To begin with,  we use again the change of variable $r=e^t$ to convert  system \eqref{vc system} of $v_c$ on $\cK$ to the system of $w$ on $\cC$, where we denote 
\[
w(t,\tht)=v_c(r,\tht).
\]
Direct computations lead to the system for $w$ as below
\beq\label{w system 1}
\cU(e^t,\p_t)w=\big(e^{2t}h_c(e^t,\tht), \,f_c(e^t),\, e^t g_c(e^t)\big),
\eeq
where the operator
\[
\cU(e^t,\p_t)=\big(\na\cdot P_w\na,\, \cdot|_{\Ga_t},\, \p^{P_w}_{n_b}\cdot|_{\Ga_b}\big).
\]
Here the coefficient matrix reads
\[
P_w(t,\tht)=P^t_\tht P_c(e^t,\tht)P_\tht\quad\hbox{with}\ P_\tht=\left(\begin{matrix}\cos\tht &-\sin \tht\\ \sin\tht& \cos\tht\end{matrix}\right)
\]
and notice that 
\[
\na=\na_{t,\tht}
\]
 in the strip domain $\cC$.

Applying Lemma \ref{space V and W} and recalling the assumptions of this proposition, we have $w\in W^2_{2, \be-1}(\cC)$ and the right side of \eqref{w system 1} satisfies
\beq\label{right side 1}
\cU(e^t,\p_t)w\in W^1_{2, \be-1}(\cC)\times W^{5/2}_{2, \be-1}(\Ga_t)\times W^{3/2}_{2, \be-1}(\Ga_b).
\eeq

\medskip
We want to show that $w\in W^3_{2, \be-1}(\cC)$, which can be done in two steps.\\
\noindent{\bf Step 1}. Localization in $t$ and standard elliptic estimates for $\zeta_k w$. Similarly as in \cite{KMR}, let $\{\zeta_k\}_{k\in\Z}\subset \cC^\infty_0(\R)$ be a partition of unity with $\zeta_k$ supported on $(k-1,k+1)$  and satisfying
\[
|\zeta^{(j)}_k(t)|< c_j,\quad\forall t\in \R,\ j\in\{0,1,2,\dots\}.
\]
Here the constant $c_j$ doesn't depend on $k,t$. Meanwhile, take
\[
\eta_k=\zeta_{k-1}+\zeta_k+\zeta_{k+1},
\]
so one has $\eta_k\zeta_k=\zeta_k$, i.e. $\eta_k=1$ on the support of $\zeta_k$.

Recalling that $w\in W^2_{2, \be-1}(\cC)$, which  implies 
\[
\zeta_k w\in H^2(\cC)
\]
 satisfying the system
\[
\cU(t,\p_t)\zeta_k w=\zeta_k\,\cU(t,\p_t)w+[\cU(t,\p_t),\,\zeta_k]w,
\]
or equivalently
\beq\label{zeta k w system}
\left\{\begin{array}{ll}
\na\cdot P_w\na(\zeta_k w)=\zeta_k e^{2t}h_c(e^t,\tht)+[\na\cdot P_w\na,\,\zeta_k]w,\qquad\hbox{on}\quad \cC\\
\zeta_kw|_{\Ga_t}=\zeta_kf_c(e^t),\quad \p^{P_w}_{n_b}(\zeta_kw)\big|_{\Ga_b}=\zeta_k  e^t g_c(e^t)+[\p^{P_w}_{n_b},\,\zeta_k]w\big|_{\Ga_b}.
\end{array}\right.
\eeq
To estimate the right side of the system above, one knows firstly from \eqref{right side 1} that
\[
\zeta_k\, \cU(t,\p_t)w\in H^1(\cC)\times H^{5/2}(\Ga_t)\times H^{3/2}(\Ga_b).
\]
On the other hand, a direct computation from \eqref{w system 1} shows 
\[
[\na\cdot P_w\na,\,\zeta_k]w=[\na\cdot P_w\na,\,\zeta_k]\eta_k\,w=\na\cdot P_w\left(\begin{array} {ll}\zeta'_k \\0\end{array}\right)(\eta_k w)+\left(\begin{array} {ll}\zeta'_k \\0\end{array}\right)\cdot P_w\na(\eta_k w)
\]
and 
\[
[\p^{P_w}_{n_b},\,\zeta_k]w\big|_{\Ga_b}={\bf n}_b\cdot P_w\left(\begin{array} {ll}\zeta'_k \\0\end{array}\right)(\eta_k w)\Big|_{\tht=-\om_2}.
\]
Consequently, one has
\[
[\cU(t,\p_t),\,\zeta_k]w=[\cU(t,\p_t),\,\zeta_k]\eta_k w\in H^1(\cC)\times H^{5/2}(\Ga_t)\times H^{3/2}(\Ga_b)
\]
satisfying  the estimate
\[\begin{split}
&\|[\cU(t,\p_t),\,\zeta_k]w\|_{H^1(\cC)\times H^{5/2}(\Ga_t)\times H^{3/2}(\Ga_b)}\\
&\le \left\|\na\cdot P_w\left(\begin{array} {ll}\zeta'_k \\0\end{array}\right)(\eta_k w)\right\|_{H^1(\cC)}
+\left\|\left(\begin{array} {ll}\zeta'_k \\0\end{array}\right)\cdot P_w\na(\eta_k w)\right\|_{H^1(\cC)}
+\left\|{\bf n}_b\cdot P_w\left(\begin{array} {ll}\zeta'_k \\0\end{array}\right)(\eta_k w)\right\|_{H^{3/2}(\Ga_b)}\\
&\le C(\|\eta'\|_{W^{2,\infty}})\|\eta_k w\|_{H^2(\cC)}
\end{split}
\]  
where Lemma \ref{restriction trace} is applied on the boundary.
Besides, one notices  that  $e^t=r$ appears in  $\na P_c(e^t,\tht)$, which can be bounded by $\del$. 

As a result, summing up the estimates above and  applying standard elliptic theories (for example Theorem 2.9 \cite{Lannes}) leads to $\zeta_kw\in H^3(\cC)$ with the estimate
\beq\label{zeta k w H3 estimate}
\begin{split}
\|\zeta_k w\|_{H^3(\cC)}&\le C(\|P_w\|_{W^{2,\infty}})\Big(\|\zeta_k\, \cU(t,\p_t)w\|_{ H^1(\cC)\times H^{5/2}(\Ga_t)\times H^{3/2}(\Ga_b)}\\
&\qquad\qquad\qquad\qquad+\|[\cU(t,\p_t),\,\zeta_k]w\|_{H^1(\cC)\times H^{5/2}(\Ga_t)\times H^{3/2}(\Ga_b)}\Big)\\
&\le C(\|\eta'\|_{W^{2,\infty}})\Big(\|\zeta_k\, \cU(t,\p_t)w\|_{ H^1(\cC)\times H^{5/2}(\Ga_t)\times H^{3/2}(\Ga_b)}+\|\eta_kw\|_{H^2(\cC)}\Big).
\end{split}
\eeq
Notice that the coefficient $C(\|\eta'\|_{W^{2,\infty}})$ above doesn't depend on $k$, which is the key to go back to the weighted norm for $w$.

\medskip
\noindent{\bf Step 2}. The estimate for $w$.  To begin with,  we will  convert the estimate above for $\zeta_k\,w$ to the estimate for $w$. In fact, one has for each $k\in \Z$ that
\[
\zeta_k w\in W^3_{2, \be-1}(\cC)
\]
from the definition of $W^3_{2,-1}(\cC)$ and $\zeta_k$. 
Moreover, it's straightforward to see that
\[
c_1\|\zeta_k \,w\|_{H^3(\cC)}\le e^{(1-\be)k}\|\zeta_k\,w\|_{W^3_{2, \be-1}(\cC)}\le c_2\|\zeta_k \,w\|_{H^3(\cC)}
\]
where $c_1, c_2$ are two constants independent of $k$. 

Consequently, multiplying $e^{(\be-1)k}$ on both sides of \eqref{zeta k w H3 estimate}, one derives 
\[
\|\zeta_k\,w\|_{W^3_{2, \be-1}(\cC)}\le C(\|\eta'\|_{W^{2,\infty}})\Big(\|\zeta_k\,\cU(t,\p_t)w\|_{W^1_{2, \be-1}(\cC)\times W^{5/2}_{2, \be-1}(\Ga_t)\times W^{3/2}_{2,-1}(\Ga_b)}+\|\eta_k\,w\|_{W^2_{2, \be-1}(\cC)}\Big).
\]
The following lemma from \cite{KMR} tells us  the relationship between  the norms of $w$ and $\zeta_k w$.
\bthm{Lemma}\label{zeta k w}
{\it Let $\{\zeta_k\}$ be the partition of unity in $\R$ defined above and $\be_0\in \R$. Then there exist positive real constants $c_1,c_2$ depending only on $l\ge 1$ such that
\[
c_1\|w\|_{W^l_{2,\be_0}(\cC)}\le \big(\sum^{+\infty}_{k=-\infty}\|\zeta_k w\|^2_{W^l_{2,\be_0}(\cC)}\big)^{1/2}\le c_2\|w\|_{W^l_{2,\be_0}(\cC)}
\]
for each $w\in W^l_{2,\be_0}(\cC)$. 
}
\ethm
Consequently, applying this lemma, we have immediately $w\in W^3_{2, \be-1}(\cC)$ with the estimate 
\[
\|w\|_{W^3_{2, \be-1}(\cC)}\le C(\|\eta'\|_{W^{2,\infty}})\Big(\|\cU(t,\p_t)w\|_{W^1_{2, \be-1}(\cC)\times W^{5/2}_{2, \be-1}(\Ga_t)\times W^{3/2}_{2, \be-1}(\Ga_b)}+\|w\|_{W^2_{2, \be-1}(\cC)}\Big).
\]
Combining Lemma \ref{space V and W} and Proposition \ref{weighted estimate near Xc}, we finish the proof for the case $l=3$.

\noindent{\bf Step 3}. The case $l>3$.  In fact, applying Theorem 2.9 \cite{Lannes} to $\zeta_k w$ system \eqref{zeta k w system}, one obtains
\[
\begin{split}
\|\zeta_k w\|_{H^l(\cC)}&\le C(\|P_w\|_{W^{l-1,\infty}})\Big(\|\zeta_k\, \cU(t,\p_t)w\|_{ H^{l-2}(\cC)\times H^{l-1/2}(\Ga_t)\times H^{l-3/2}(\Ga_b)}\\
&\qquad\qquad\qquad\qquad+\|[\cU(t,\p_t),\,\zeta_k]w\|_{H^1(\cC)\times H^{5/2}(\Ga_t)\times H^{3/2}(\Ga_b)}\Big)\\
&\le C(\|\eta'\|_{W^{l-1,\infty}})\Big(\|\zeta_k\, \cU(t,\p_t)w\|_{ H^{l-2}(\cC)\times H^{l-1/2}(\Ga_t)\times H^{l-3/2}(\Ga_b)}+\|\eta_kw\|_{H^2(\cC)}\Big).
\end{split}
\]
The rest part can be proved similarly as well and the proof is finished.
\ef


\subsection{Proof of Theorem \ref{main thm}}
Now we are ready to prove this main theorem.  First of all, recalling definition \eqref{def of V} of the weighted space $V^l_\be(\Om)$, one knows that $u$ is divided into  $v_c$ and $v_R$.  Therefore,  the proof deals with these two parts and an inductive method is applied here for $l\ge 2$.

\medskip
\noindent{\bf Step 1: $l=2$.} For the key part concerning $v_c$,   we  apply Proposition \ref{order-l regularity} directly if the assumptions there are satisfied. In fact, checking from system \eqref{vc system}, one can see that 
\[
h_c=\big(\chi_c \,h-[\chi_c,\Del]u\big)\circ T_c\in V^{0}_{\be}(\cK).
\]
This holds since one has 
\[
\chi_c h\circ T_c\in V^{0}_{\be}(\cK)
\]
by the assumption $h\in V^{0}_{\be}(\Om)$ and moreover $u\in H^2(\Om)$ leads to 
\[
r^\be[\chi_c,\Del]u\circ T_c \in L^2(\cK)
\]
from  definition \eqref{def of V}. 

On the other hand, $f\in V^{3/2}_\be(\Ga_t)$ implies  $f_c\in V^{3/2}_\be(\Ga_t)$ immediately, so it remains to check $g_c$. Recalling that
\[
g_c=\Big(\chi_c|_{\Ga_b}g-(\p_{n_b}\chi_c)|_{\Ga_b}u\Big)\circ T_c,
\]
one obtains $\chi_c|_{\Ga_b}g\circ T_c\in V^{1/2}_\be(\Ga_b)$ directly from the assumption of this theorem. Meanwhile, one also has $u\in H^{3/2}(\Ga_b)$, which infers 
\[
(\p_{n_b}\chi_c)|_{\Ga_b}u\circ T_c\in H^{1/2}(\Ga_b).
\]
Consequently, checking \eqref{boundary norm} for the norm of $V^{1/2}_\be(\Ga_b)$ and noticing that 
\[
(\p_{n_b}\chi_c)|_{\Ga_b}u\circ T_c=0\quad \hbox{near}\ X_c
\]  one derives  $(\p_{n_b}\chi_c)|_{\Ga_b}u\circ T_c\in V^{1/2}_\be(\Ga_b)$ immediately. 

In a word, the assumptions of Proposition \ref{weighted estimate near Xc} are all satisfied indeed.  
Applying this proposition, we have  from estimate \eqref{vc V20 estimate}  that
\[
\begin{split}
\|v_c\|_{v^2_{\be}(\cK)}&\le C\big(\|h_c\|_{V^{0}_{\be}(\cK)}+\|f_c\|_{V^{3/2}_{\be}(\Ga_t)}+\|g_c\|_{V^{1/2}_{\be}(\Ga_b)}\big)\\
&\le C\Big(\|h\|_{V^0_\be(\cK)}+|f|_{V^{3/2}_\be(\Ga_t)}+|g|_{V^{1/2}_\be(\Ga_b)}
+\del\|u\|_{V^2_\be(\cK)}\Big),
\end{split}
\]
where the constant $C$ depends on $\cK,\chi_c$,  and the following inequalities are applied:
\[
\begin{split}
\|[\chi_c,\Del]u\circ T_c\|_{V^0_\be(\cK)}&\le C\big(\|r^\be(\Del\chi_c)\bar\chi_c u\circ T_c\|_{L^2(\cK)}+\|r^\be\na\chi_c\cdot\na(\bar \chi_c u)\circ T_c\|_{L^2(\cK)}\big)\\
&\le \del\,C\|u\|_{V^2_\be(\cK)}
\end{split}
\]
and 
\[
\|(\p_{n_b}\chi_c)|_{\Ga_b}u\circ T_c\|_{V^{1/2}_\be(\Ga_b)}\le C\|(\p_{n_b}\chi_c)\bar\chi_cu\circ T_c\|_{V^1_\be(\cK)}
\le \del\, C\|u\|_{V^2_\be(\cK)}
\]
with $\bar \chi_c$ another $C^\infty_0(\Om)$ function satisfying $\chi_c=\bar \chi_c\chi_c$.

Secondly, for the remainder $v_R=(1-\chi_c)u\circ T_R$, direct computations show that $v_R$ satisfies the following system 
\beq\label{vR system}
\left\{\begin{array}{ll}\na\cdot P_R\na v_R=h_R\qquad \hbox{on}\quad R\\
v_R|_{\Ga_t}= f_R,\quad \p^{P_R}_{n_b}v_R|_{\Ga_b}= g_R\end{array}\right.
\eeq where
\[
h_R=\Big((1-\chi_c) \,h-[1-\chi_c,\Del]u\Big)\circ T_R,\ f_R=\big((1-\chi_c)|_{\Ga_t}\,f\big)\circ T_R,
\]
and 
\[
g_R=\big(1+l'(x)^2\big)^{1/2}\Big((1-\chi_c)|_{\Ga_b}\,g+(\p_{n_b}\chi_c)u|_{\Ga_b}\Big)\circ T_R.
\]
Moreover, the coefficient matrix reads
\[
P_R=(P_r)^tP_r\quad\hbox{with}\quad P_r=(\na T^{-1}_R)\circ T_R=\left(\begin{matrix} 1& -\big(\eta(x)-l(x)\big)^{-1}\big((\eta'(x)-l'(x))z+l'(x)\big)\\ 0 & \big(\eta(x)-l(x)\big)^{-1}\end{matrix}\right).
\] 
Since this system for $v_R$ is defined on the flat strip $R$, standard elliptic theories apply directly (for example Theorem 2.9 \cite{Lannes}). Meanwhile, one has immediately that $v_R\in H^2(R)$ with the estimate 
\[
\begin{split}
\|v_R\|_{H^2(R)}&\le C(\|\eta'\|_{W^{1,\infty}})\big(\|h_R\|_{L^2(R)}+\|f_R\|_{H^{3/2}(\Ga_t)}+\|g_R\|_{H^{1/2}(\Ga_b)}\big)\\
&\le C(\|\eta'\|_{W^{1,\infty}})\Big(\|h\circ T_R\|_{L^2(R)}+\|f\circ T_R\|_{H^{3/2}(R)}+\|g\circ T_R\|_{H^{1/2}(R)}
+\|u\|_{H^1(\Om)}\Big)
\end{split}
\]
and notice that the term $\|u\|_{H^1(\Om)}$ can be handled using $H^2$ estimate from \cite{MW1} and Lemma \ref{H32 V32_0}:
\[
\begin{split}
\|u\|_{H^1(\Om)}&\le C(\|\eta'\|_{W^{1,\infty}})\big(\|h\|_{L^2(\Om)}+\|f\|_{H^{3/2}(\Ga_t)}+\|g\|_{H^{1/2}(\Ga_b)}\big)\\
&\le C(\|\eta'\|_{W^{1,\infty}})\Big(\|h\|_{V^0_0(\Om)}+|f|_{V^{3/2}_0(\Ga_t)}+|g|_{V^{1/2}_0(\Ga_b)}\Big).
\end{split}
\]

As a result, combining these estimates above, we have proved that $u\in V^2_\be(\Om)$ satisfies the estimate 
\beq\label{u V20 estimate}
\|u\|_{V^2_\be(\Om)}\le C(\|\eta'\|_{W^{1,\infty}})\Big(\|h\|_{V^0_\be(\Om)}+|f|_{V^{3/2}_\be(\Ga_t)}+|g|_{V^{1/2}_\be(\Ga_b)}\Big).
\eeq

\medskip
\noindent{\bf Step 2:} Case (i) $l\ge 3$ when $\eta\in W^{l,\infty}(\R^+)$.  An induction method is used in this part.  To begin with, we know already from the assumptions of this theorem that
\[
h\in V^{l-2}_{l-2+\be}(\Om),\ f\in V^{l-1/2}_{l-2+\be}(\Ga_t),\ g\in V^{l-3/2}_{l-2+\be}(\Ga_b).
\] 
Assuming  $u\in V^{l-1}_{l-3+\be}(\Om)$ and the following estimate holds
\[
\|u\|_{V^{l-1}_{l-3+\be}(\Om)}\le C(\|\eta'\|_{W^{l-2,\infty}})\big(\|h\|_{V^{l-3}_{l-3+\be}(\Om)}+|f|_{V^{l-3/2}_{l-3+\be}(\Ga_t)}+|g|_{V^{l-5/2}_{l-3+\be}(\Ga_b)}\big),
\]
 we are going to show that $u\in V^l_{l-2+\be}(\Om)$ with  corresponding estimate.

Firstly, we deal with the part $v_c$. In fact,  it remains again to check the assumptions of Proposition \ref{order-l regularity}.  Since one has $u\in V^{l-1}_{l-3+\be}(\Om)\cap H^2(\Om)$, one deduces directly that
\[
[\chi_c,\Del]u\circ T_c\in V^{l-2}_{l-2+\be}(\cK)
\] 
and direct computations lead to
\[
\begin{split}
\|[\chi_c,\Del]u\circ T_c\|_{V^{l-2}_{l-2+\be}(\cK)}&\le C\Big(\|(\Del \chi_c)\bar\chi_cu\circ T_c\|_{V^{l-2}_{l-2+\be}(\cK)}+\|\big(\na \chi_c\cdot\na (\bar \chi_cu)\big)\circ T_c\|_{V^{l-2}_{l-2+\be}(\cK)}\Big)\\
&\le C(\|\eta'\|_{W^{l-2,\infty}}) \|u\|_{V^{l-1}_{l-3+\be}(\Om)},
\end{split}
\]
where one notices that $\Del \chi_c, \na\chi_c$  vanish near $X_c$ and $C(\|\eta'\|_{W^{l-2,\infty}})$ comes from $\na (\bar \chi_cu)\circ T_c$.
Consequently, one has $h_c\in V^{l-2}_{l-2+\be}(\cK)$ satisfying
\[
\|h_c\|_{V^{l-2}_{l-2+\be}(\cK)}\le C(\|\eta'\|_{W^{l-2,\infty}})\big(\|h\|_{V^{l-2}_{l-2+\be}(\Om)}+\|u\|_{V^{l-1}_{l-3+\be}(\Om)}\big).
\]

On the other hand, the definition of $V^{l-1/2}_{l-2+\be}(\Ga_t)$ infers immediately that $f_c\in V^{l-1/2}_{l-2+\be}(\Ga_t)$ holds. Meanwhile, for the term $g_c$, one can show directly from Lemma \ref{restriction trace} and Lemma \ref{imedding} that 
\[
\begin{split}
\|g_c\|_{V^{l-3/2}_{l-2+\be}(\Ga_b)}&
\le C\Big(\|g\|_{V^{l-3/2}_{l-2+\be}(\Ga_b)}+\|(\p_{n_b}\chi_c)|_{\Ga_b} u\circ T_c\|_{V^{l-3/2}_{l-2+\be}(\Ga_b)}\Big)
\\
&\le C\Big(\|g\|_{V^{l-3/2}_{l-2+\be}(\Ga_b)}+\|\bar \chi_c u\circ T_c\|_{V^{l-1}_{l-2+\be}(\cK)}\Big)\le C\Big(\|g\|_{V^{l-3/2}_{l-2+\be}(\Ga_b)}+\|u\|_{V^{l-1}_{l-3+\be}(\Om)}\Big).
\end{split}
\]
 Summing these up,  one can see that the assumptions of Proposition \ref{order-l regularity} are all satisfied. Applying Proposition \ref{order-l regularity} and Lemma \ref{imedding} together with \eqref{u V20 estimate}, we finally have $v_c\in V^l_{l-2+\be}(\cK)$ with the estimate
\[
\begin{split}
\|v_c\|_{V^l_{l-2+\be}(\cK)}&\le C(\|\eta'\|_{W^{l-1,\infty}})\Big(\|h\|_{V^{l-2}_{l-2+\be}(\Om)}+\|f\|_{V^{l-1/2}_{l-2+\be}(\Ga_t)}+\|g\|_{V^{l-3/2}_{l-2+\be}(\Ga_b)}+\|u\|_{V^{l-1}_{l-3+\be}(\Om)}\Big)\\
&\le C(\|\eta'\|_{W^{l-1,\infty}})\Big(\|h\|_{V^{l-2}_{l-2+\be}(\Om)}+\|f\|_{V^{l-1/2}_{l-2+\be}(\Ga_t)}+\|g\|_{V^{l-3/2}_{l-2+\be}(\Ga_b)}\Big).
\end{split}
\]
Secondly, applying standard elliptic theories, one derives an estimate for $v_R\in H^l(R)$ similarly as before.  

As a result, combining these two parts together, one  concludes that $u\in V^l_{l-2+\be}(\Om)$ satisfying
the desired estimate.
\medskip

\noindent{\bf Step 3:} Case (ii) $l\ge 3$ when $\eta\in H^{l-1/2}(\R^+)$. Here we need the regularizing transformation $\tilde T_c$ defined in \eqref{def of tilde Tc} in Section 2.1. Therefore, we change every $T_c$ we meet  into $\tilde T_c$, and the corresponding coefficient matrix $P_c$ in system \eqref{vc system} of $v_c$ should be replaced by 
\[
\tilde P_c=(\tilde P^{-1}_0)^t(\tilde P_S\circ \tilde T_0)^t\tilde P_S\circ T_0(\tilde P_0)^{-1},
\]
where  $\tilde P_S\circ \tilde T_0$ is given  by \eqref{tilde PS}. So one can tell directly that 
\[
\tilde P_c=\tilde P_c(\na s).
\]
Similarly as before,  one can check the estimates from the beginning to find out that all we need is to focus on system \eqref{zeta k w system} for $\zeta_k w$ again.  The corresponding coefficient matrix $P_w$ should be replaced by 
\[
\tilde P_w=P^t_\tht\tilde P_c\big(\na s(e^t,\tht)\big)P_\tht=\tilde P_{w,1}+\tilde P_{w,2}
\] with
\[
\tilde P_{w,1}=P^t_\tht\tilde P_c({\bf 0})P_\tht,\quad \tilde P_{w,2}=P^t_\tht\big(\tilde P_c\big(\na s(e^t,\tht)\big)-\tilde P_c({\bf 0})\big)P_\tht.
\]
Applying Theorem 2.9 ii)\cite{Lannes} and assuming that $l\ge 3$, one has 
\beq\label{zeta k w estimate 1}\begin{split}
\|\zeta_k w\|_{H^l(\cC)}\le C\big(l, \|\tilde P_{w,1}\|_{W^{l-1,\infty}(\cC)},\|\tilde P_{w,2}\|_{H^{l-1}(\cC)}\big)&\Big(\|\zeta_k\, \cU(t,\p_t)w\|_{ H^{l-2}(\cC)\times H^{l-1/2}(\Ga_t)\times H^{l-3/2}(\Ga_b)}\\
&\quad+\|[\cU(t,\p_t),\,\zeta_k]w\|_{H^{l-2}(\cC)\times H^{l-1/2}(\Ga_t)\times H^{l-3/2}(\Ga_b)}\Big)
\end{split}
\eeq
To handle the second term on the right side of the inequality above,   the following three terms 
\[
\left\|\na\cdot \tilde P_w\left(\begin{matrix} \zeta'_k\\ 0\end{matrix}\right)w\right\|_{H^{l-2}(\cC)},\  \left\|\left(\begin{array} {ll}\zeta'_k \\0\end{array}\right)\cdot \tilde P_w\na(\eta_k w)\right\|_{H^{l-2}(\cC)}\ \hbox{and}\ 
\left\|{\bf n}_b\cdot \tilde  P_w\left(\begin{array} {ll}\zeta'_k \\0\end{array}\right)(\eta_k w)\right\|_{H^{l-3/2}(\Ga_b)}
\]
need to be taken care of according to the proof of Propostion \ref{order-l regularity}.  

For example, in the first term above, we consider the estimate for the term where all the derivatives are taken on one $\na s(e^t,\tht)$. In fact,  this  term is like
\[
A=P^t_\tht\, (\p \tilde P_c)\big(\na s(e^t\cos \tht, e^t\sin \tht)\big)\,P_\tht\,(\p^\al\p s)(e^t\cos \tht, e^t\sin \tht) e^{(l-1)t} \phi(\cos \tht,\sin \tht)\left(\begin{matrix} \zeta'_k\\ 0\end{matrix}\right)w,
\]
which comes from $\p^\al \tilde P_w$ with $|\al|=l-1$.  Then the following estimate  holds:
\[
\begin{split}
\|A\|_{L^2(\cC)}&\le C \|(\p^\al\p s)(e^t\cos \tht, e^t\sin \tht) e^{(l-1)t}\|_{L^2(\cC)}\|\eta_kw\|_{L^\infty(\cC)}\\
&\le C\|s\|_{H^l(\cK)}\|\eta_kw\|_{H^{l-1}(\cC)},
\end{split}
\]
where we notice that $e^{(l-1)t}$ is used to transform between different domains.
Similarly,  one can have estimates for the other terms in \eqref{zeta k w estimate 1}. Plugging all the estimates back into \eqref{zeta k w estimate 1}, one derives
\[
\|\zeta_k w\|_{H^l(\cC)} \le C\big(\|s\|_{H^l(\cK)}\big)\|\eta_k w\|_{H^{l-1}(\cC)}.
\]
As a result, remembering the definition of $s$ and applying \eqref{trace for s}, one finds immediately that 
\[
\|s\|_{H^l(\cK)}\le C\big(|\eta|_{H^{l-1/2}(\R^+)}\big).
\] 
 Therefore, the proof of Theorem \ref{main thm} ii) can be finished. 
\ef

\bigskip

\section{Some other boundary-value problems}
\subsection{Dirichlet boundary problem}
To consider about the system $\mbox{(DVP)}$,  one needs to verify first of all the existence of the solution $u\in H^2(\Om)$.  To begin with, it is well-known from Lemma 4.4.3.1 \cite{PG} that there exists a unique variational solution $u\in H^1(\Om)$ under the assumptions of Theorem \ref{D and N thm}.

Secondly, we recall an adjusted version of Theorem 3.2.5 \cite{BR} here.
\bthm{Theorem}\label{BR thm}(Banasiak-Roach1989) {\it Let \\
(i) $\Om_c\subset \R^2$ be a bounded, open set with only one angle and a $C^{2,0}$ curvilinear polygon as its boundary;\\
(ii) The boundary values $f,g$ satisfy the following conditions:\\
For system $\mbox{(DVP)}$,  one assumes that $f\in H^{3/2}(\Ga_t)$, $g\in H^{3/2}(\Ga_b)$ and 
\[
f|_{X_c}=g|_{X_c};
\]
For system $\mbox{(NVP)}$, one assumes that $f\in H^{1/2}(\Ga_t)$, $g\in H^{1/2}(\Ga_b)$.

Then the variational solution $u\in H^1(\Om)$ of system $\mbox{(DVP)}$ or system $\mbox{(NVP)}$ can be uniquely represented in the form
\[
u=u_r+\sum_{m\in L} c_{m}S_{m}
\]
with $u_r\in H^2(\Om_c)$ and some $c_{m}\in \R$. Moreover, the set $L$ is defined as
\[
L=\{m|\, -1<\lam_m<0\},\quad\hbox{where}\quad \lam_m=\f{m\pi}{\om},
\]
and $\om $ is the angle.
}
\ethm

Therefore, localizing system $\mbox{(DVP)}$ as in the mixed-boundary case and checking directly from this theorem, it is clear that  when the contact angle $\om\in(0,\pi)$, one has 
\[
\lam_m\notin (-1,0),
\]
which implies that there exists a unique solution $u\in H^2(\Om)$ to $\mbox{(DVP)}$.  By the way, the compatibility condition  can also be found in \cite{PG1, MW1}.

Consequently, the solution $u$ lies  in the space $V^2_2(\Om)$ as before.  To prove the first part in Theorem \ref{D and N thm}, we only needs to follow the proof in Section 3 and check line by line. In this case, we will consider the Dirichlet boundary system for $v_c$:
\[
\left\{\begin{array}{ll}\na\cdot P_c\na v_c=h_c\qquad \hbox{on}\quad \cK\\
v_c|_{\Ga_t}= f_c,\quad v_c|_{\Ga_b}=g_c,\end{array}\right.
\]
and the system for $w(t,\tht)=v_c(e^t,\tht)$ is  slightly different as well.

In Proposition \ref{singularity decomposition}, one also has a different eigenvalue problem for $\cU(\lam)$ here:
\[
\left\{\begin{array}{ll}
-\phi''(\tht)-\lam^2\phi(\tht)=0,\quad \tht\in I\\
\phi(\om_1)=0,\qquad \phi(-\om_2)=0
\end{array}\right.
\]
with the eigenvalues and eigenfucntions 
\[
\lam_m=\f{m\pi}{\om},\quad \phi_m(\tht)=\sin{\big(\lam_m(\tht+\om_2)\big)}\qquad\hbox{for}\ \forall m\in\Z.
\]
Consequently, we can tell that the eigenvalues 
\[
\lam_m\notin [-1,0)\quad \hbox{for}\, \om\in (0,\pi), \forall m\in\Z.
\]
Therefore, following the proof of Proposition \ref{singularity decomposition} and Proposition \ref{weighted estimate near Xc}, we  conclude that  
\[
v_c\in V^2_{1+\be}(\cK)\quad\hbox{ for any }\be\in(0,1],
\]
 since we cannot cross over $0$ for the eigenvalue $\lam_m$. 

Moreover,  the regularity of $v_c$ can be improved in the same way as before, and the weighted estimates rely on standard elliptic estimates with Dirichlet boundaries (which can be proved similarly as in \cite{Lannes}). 

As a result, the proof for the Dirichlet case in Theorem \ref{D and N thm} is finished.

\subsection{Neumann boundary problem} The case of Neumann boundaries can be proved similarly , and the eigenvalue value $\lam_m$ from the eigenvalue problem  for $\cU(\lam)$ turns out to be the same as  in the Dirichlet case. 

Here, one needs to notice that we assume the existence of the solution $u$ when $\Om$ is unbounded.  When the domain is bounded, the existence can be proved under the compatibility condition for $\mbox{(NVP)}$.

Therefore, the proof for the second part of Theorem \ref{D and N thm} follows.

\section{Application on the Dirchlet-Neumann operator}
In the end, we  show that  the weighted elliptic theory above can be applied to the Dirichlet-Neumann operator, which is  an important operator in the water-waves problem.  

To begin with, recalling that for a proper function $f$ on $\Ga_t$, the D-N operator $\cN$ is defined as
\beq\label{def of DN op}
\cN f:=\na_{n_t}f_{\cH}\big|_{\Ga_t}
\eeq 
where $f_{\cH}$ is the harmonic extension of $f$ satisfying the system
\[
\left\{\begin{array}{ll}\Del f_{\cH}=0,\qquad \hbox{on}\quad \Om,\\
f_{\cH}|_{\Ga_t}=f,\qquad \na_{n_b}f_{\cH}|_{\Ga_b}=0.
\end{array}
\right.
\]

When the function $f\in V^l_{l-2+\be}(\Ga_t)$ for any $\be\in[0,2]$, we know directly from Theorem \ref{main thm} that 
\beq\label{f_H estimate}
\|f_{\cH}\|_{V^l_{l-2+\be}(\Om)}\le C(\|\eta\|_{H^{l-1/2}(\R^+)})\|f\|_{V^{l-1/2}_{l-2+\be}(\Ga_t)}.
\eeq
On the other hand, to consider the estimate for the D-N operator, we need the following lemma about the product estimate in the weighted space :
\bthm{Lemma}\label{product estimate}{\it For any two functions $f\in V^{k+1/2}_\be(\Ga_t)$ and $g\in H^{k+1/2}(\Ga_t)$ with an integer $k\ge 2$ and a real $\be$, one has  the  estimate for the product of $f,\,g$:
\[
\|f\,g\|_{V^{k+1/2}_\be(\Ga_t)}\le C \|f\|_{V^{k+1/2}_\be(\Ga_t)}\|g\|_{H^{k+1/2}(\Ga_t)}
\]
where the constant $C$ depends only on $k$.
}
\ethm

\noindent{Proof}. According to the definition of $V^{k+1/2}_\be(\Ga_t)$, it suffices to show that the estimate holds for $f_c,\,g_c$ defined on $\Ga_t$ of $\cK$ near the contact point. 

In fact, checking  \eqref{boundary norm}, one knows immediately that
\[\begin{split}
\|f_c\,g_c\|_{V^{k+1/2}_\be(\Ga_t)}\le &C(k)\sum_{j\le k}\int_{R^+}r^{2(\be-k)-1+2j}\big|\p^j_r\big(f_c(r)g_c(r)\big)\big|^2dr\\
&
+C(k)\sum_{j\le k}\int_{\R^+}\int^{2r}_{r/2}r^{2(\be-k)}\f{\big|r^j\p^j_r(f_cg_c)(r)-\rho^j\p^j_\rho(f_cg_c)(\rho)\big|^2}{|r-\rho|^2}d\rho dr\\
:=&C\sum _{j\le k} (A_j+B_j)
\end{split}
\]
where $A_j$, $B_j$ will be handled one by one.

Firstly, one has for the term $A_j$ that
\beq\label{A_j}
A_j\le \sum_{0\le\al\le j}C_\al\int_{R^+}r^{2(\be-k)-1+2j}\big|\p^\al_rf_c(r)\,\p^{j-\al}_rg_c(r)\big)\big|^2dr,
\eeq
which is separated into two cases. 

For the first case when $j-\al=k$, i.e. $j=k\ge  2,\,\al=0$,  all the derivatives are taken on $g$, so the corresponding term becomes
\[\begin{split}
\int_{R^+}r^{2(\be-k)-1+2k}\big|f_c(r)\,\p^k_rg_c(r)\big)\big|^2dr
&\le C(\del)\int_{R^+}r^{2(\be-k)-1}|rf_c(r)|^2|\p^k_rg_c(r)|^2 dr\\
&\le C \|r^{\be-k-1/2} r f_c\|^2_{L^\infty(R^+)}\|\p^k_rg_c\|^2_{L^2(\R^+)}\\
&\le C\|r^{\be-k-1/2} r f_c\|^2_{H^1(\R^+)}\|g_c\|^2_{H^k(\R^+)}\le C\|f_c\|^2_{V^{k+1/2}_\be(\Ga_t)}\|g_c\|^2_{H^k(\R^+)}
\end{split}\]
where  the imbedding theorem from $H^1(\R^+)$ to $L^\infty(\R^+)$ is applied.

For the other case when $j-\al\le k-1$, the corresponding terms are  handled  in a slightly different way:
\[\begin{split}
\int_{R^+}r^{2(\be-k)-1+2j}\big|\p^\al_rf_c(r)\,\p^{j-\al}_rg_c(r)\big)\big|^2dr
&\le \int_{R^+}r^{2(\be-k)-1+2j}\big|\p^\al_rf_c(r)\big|^2dr \,\|\p^{j-\al}_rg_c\|_{L^\infty(\R^+)}\\
&\le C(\del) \|f_c\|^2_{V^{k+1/2}_\be(\Ga_t)}\|g_c\|^2_{H^k(\R^+)},
\end{split}\]
where  the imbedding theorem from $H^1(\R^+)$ to $L^\infty(\R^+)$ is applied again.
Consequently, substituting the two cases above into \eqref{A_j}, one arrives at 
\beq\label{A_j sum}
A_j\le C(\del)\,\|f_c\|^2_{V^{k+1/2}_\be(\Ga_t)}\|g_c\|^2_{H^k(\R^+)}.
\eeq

Secondly, to deal with $B_j$, one deduces that
\beq\label{B_j}\begin{split}
B_j&\le \sum_{\al\le j}C_\al\int_{R^+}\int^{2r}_{r/2}r^{2(\be-k)}\f{\big|r^\al\p^\al_rf_c(r)r^{j-\al}\p^{j-\al}_rg_c(r)-\rho^\al\p^\al_\rho f_c(\rho)\rho^{j-\al}\p^{j-\al}_\rho g_c(\rho)\big|^2}{|r-\rho|^2}d\rho dr\\
&\le \sum_{\al\le j}C_\al\int^\del_{0}\int^{2r}_{r/2}r^{2(\be-k)}
\f{\big|r^{\al}\p^\al_rf_c(r)-\rho^\al \p^\al_\rho f_c(\rho)\big|^2}{|r-\rho|^2}|r^{j-\al}\p^{j-\al}_rg_c(r)|^2\\
&\qquad+\sum_{\al\le j}C_\al\int^{2\del}_{0}\int^{2r}_{r/2}r^{2(\be-k)}|\rho^\al \p^\al_\rho f_c(\rho)|^2\f{\big|r^{j-\al}\p^{j-\al}_rg_c(r)-\rho^{j-\al}\p^{j-\al}_\rho g_c(\rho)\big|^2}{|r-\rho|^2}d\rho dr
\end{split}
\eeq
where all the terms  can be treated similarly as in the case of $A_j$. 

For example, when $j=k\ge 2,\al=0$, there is a term from \eqref{B_j} satisfying
\[\begin{split}
\int^\del_{0}\int^{2r}_{r/2}r^{2(\be-k)}|\p_rf_c(t)|^2|r^k\p^k_rg_c(r)|^2 d\rho dr
&\le C(\del)\|t^{\be-k-1/2}t^2\p_rf_c(t)\|^2_{L^\infty(R^+)}\|g_c\|^2_{H^r(\R^+)}\\
&\le C\|f_c\|^2_{V^{k+1/2}_\be(\Ga_t)}\|g_c\|^2_{H^k(\R^+)},
\end{split}\]
where $t$ is some number between $r$ and $\rho$ in the integral, and the imbedding theorem from $H^1(\R^+)$ to $L^\infty(\R^+)$ is applied on $f_c$ one more time. We omit the estimates for the remainder terms since one only needs to check from one term to another  similarly as before. As a result, The proof is finished.
\ef

\bigskip
Now we conclide the weighted estimate for the D-N operator.
\bthm{Proposition}\label{D-N op}{\it 
Let $k\ge 2$ be an integer and $\be\in[k,k+2]$ be real. For any function $f\in V^{k+3/2}_\be(\Ga_t)$, one has the following estimate for the D-N operator $\cN$:
\[
\|\cN f\|_{V^{k+1/2}_\be(\Ga_t)}\le C\big(\|\eta\|_{H^{k+3/2}(\R^+)}\big)\|f\|_{V^{k+3/2}_\be(\Ga_t)}.
\]
}
\ethm

\noindent{Proof}.  The proof is a direct application of Theorem \ref{main thm} and  Lemma \ref{product estimate}.
\ef

\end{document}